\documentclass[journal]{IEEEtran}
\usepackage{subfigure}
\usepackage{algorithm}
\usepackage{algorithmic}
\usepackage{color}
\usepackage{tabularx}
\usepackage{amsmath}
\usepackage{amssymb}
\usepackage{bm}
\usepackage{graphicx}
\usepackage{booktabs}  
\usepackage{threeparttable} \usepackage{multirow}
\usepackage{cite}
\usepackage[square, comma, sort&compress, numbers]{natbib}

\ifCLASSINFOpdf
\else
\fi
\begin{document}
\title{Federated Distributionally Robust Optimization with
Non-Convex Objectives: Algorithm and Analysis}
\author{Yang Jiao, Kai Yang,~\IEEEmembership{Senior Member,~IEEE,}
        Dongjin Song,~\IEEEmembership{Member,~IEEE}


 }
 
\maketitle
\begin{abstract}
Distributionally Robust Optimization (DRO), which aims to find an optimal decision that minimizes the worst case cost over the ambiguity set of probability distribution, has been widely applied in diverse applications, \textit{e.g.}, network behavior analysis, risk management, \textit{etc.} Nevertheless, prevailing DRO techniques encounter three primary challenges in distributed environments: 1) addressing asynchronous updating efficiently; 2) leveraging the prior distribution effectively; 3) appropriately adjusting the degree of robustness based on varying scenarios. To this end, we propose an asynchronous distributed algorithm, named \textbf{A}synchronous \textbf{S}ingle-loo\textbf{P} alternat\textbf{I}ve g\textbf{R}adient proj\textbf{E}ction (ASPIRE) algorithm with the it\textbf{E}rative \textbf{A}ctive \textbf{S}\textbf{E}t method (EASE) to tackle the federated distributionally robust optimization (FDRO) problem. In addition, a new uncertainty set, \textit{i.e.}, constrained $D$-norm uncertainty set, is developed to effectively leverage the prior distribution and flexibly control the degree of robustness. We further enhance the proposed framework by integrating various uncertainty sets and conducting a comprehensive theoretical analysis of the computational complexity associated with each uncertainty set. To expedite convergence speed, we also introduce ASPIRE-ADP, a method that can dynamically adjust the number of active workers.  Finally, our theoretical analysis elucidates that the proposed algorithm is guaranteed to converge and the iteration complexity and communication complexity are also analyzed. Extensive empirical studies on real-world datasets validate that the proposed method excels not only in achieving fast convergence and robustness against data heterogeneity and malicious attacks but also in effectively managing the trade-off between robustness and performance. 
\end{abstract}
 
\begin{IEEEkeywords}
 Distributionally Robust Optimization, Distributed Optimization, Federated Learning, Uncertainty Set, Complexity Analysis
\end{IEEEkeywords}
\IEEEpeerreviewmaketitle
 
\section{Introduction}

\IEEEPARstart{O}{ver} the past decade, there has been a notable surge in the ubiquity of smartphones and Internet of Things (IoT) devices, resulting in the substantial generation of voluminous data on a daily basis. The prevailing paradigm of centralized machine learning necessitates the aggregation of data onto specific servers for the purpose of training models, thereby incurring high communication overhead \citep{sun2019communication} and suffering privacy risks \citep{sicari2015security}. As a remedial approach, distributed machine learning methods, exemplified by federated learning \cite{mohri2019agnostic,han2024federated,dubey2025understanding,cao2024overview,cao2021optimized,gupta2024unleashing,gupta2025privacy}, have been proposed. Considering a distributed system composed of $N$ workers (devices), we denote the dataset of these workers as $\{ {D_1}, \cdots,{D_N}\} $. For the $j^{\rm{th}}$ ($1 \! \le \! j \!\le \! N$) worker, the labeled dataset is given as ${D_j} = \{ {\bf{x}}^i_j,y^i_j\} $, where ${\bf{x}}^i_j \! \in \! \mathbb{R}^{d}$ and $y^i_j \in  \{1,\cdots,c\}$ denote the $i^{\rm{th}}$ data sample and the corresponding label, respectively. The distributed machine learning tasks can be formulated as the following optimization problem,
\begin{align}
\label{eq:0}
\mathop {\min }\limits_{\boldsymbol{w}\in {\boldsymbol{\mathcal{W}}}} \; F(\boldsymbol{w}) \quad {\rm{with}} \quad F(\boldsymbol{w}): = \sum\nolimits_{j} {{f_j}} (\boldsymbol{w}),
\end{align}
where $\boldsymbol{w}\in {\mathbb{R}^p}$ is the model parameter to be learned and ${\boldsymbol{\mathcal{W}}}\!\subseteq\! \mathbb{R}^p$ is a nonempty closed convex set, ${f_j}( \cdot )$ is the empirical risk over the $j^{\rm{th}}$ worker involving only the local data:
\begin{equation}
\label{eq:2}
\begin{aligned}
{f_j}({\boldsymbol{w}}) = \sum\nolimits_{i:{\bf{x}}_j^i \in {D_j}}^{} {\frac{1}{{|{D_j}|}}} \mathcal{L}_j({\bf{x}}_j^i,y_j^i;{\boldsymbol{w}}),
\end{aligned}
\end{equation}
where $\mathcal{L}_j$ is the local objective function over the $j^{\rm{th}}$ worker. Problem in Eq. (\ref{eq:0}) arises in numerous areas, such as federated learning \cite{mohri2019agnostic}, Internet of Things \cite{dubey2025integrating}, distributed signal processing \citep{geraci2015energy}, multi-agent optimization \citep{nedic2009distributed,ning2023multi}, \textit{etc}. Nevertheless, such problems fail to account for the data heterogeneity \citep{qian2020robustness}  among different workers (\textit{i.e.}, data distribution of workers could be substantially different from each other \citep{singhal2021federated}). \textcolor{black}{Studying data heterogeneity in federated learning is critically important. For example, in real-world medical imaging, hospitals often collect chest X-rays under different conditions, leading to Non-IID data distributions across institutions \cite{haripriya2025privacy}. Such heterogeneity can significantly degrade model performance if not properly addressed.} Indeed, it has been shown that traditional federated approaches, such as FedAvg \citep{mcmahan2017communication}, built for independent and identically distributed (IID) data may perform poorly when applied to Non-IID data \citep{karimireddy2019scaffold}. This issue can be mitigated via learning a robust model with the primary objective of attaining consistently high performance across all workers by solving the following distributionally robust optimization (DRO) problem in a distributed manner:
\begin{align}
\label{eq:4_new}
\mathop {{\rm{min}}}\limits_{\boldsymbol{w}\in {\boldsymbol{\mathcal{W}}}} {\rm{ }}\mathop {{\rm{max}}}\limits_{{\bf{p}} \in {\bf{\Omega} } \subseteq {\Delta _N}} F(\boldsymbol{w},{\bf{p}}): = \sum\nolimits_{j} p_j {{f_j}} (\boldsymbol{w}),
\end{align}
where ${\bf{p}}  =  [{p_1}, \cdots ,{p_N}] \! \in \! {\mathbb{R}^N}$ is the adversarial distribution in $N$ workers, the $j^{\rm{th}}$ entry in this vector, \textit{i.e.}, ${p_j}$ represents the adversarial distribution value for the $j^{\rm{th}}$ worker. ${\Delta _N}=\{ {\bf{p}} \in {\mathbb{R}^{N}_{+}}: {\bf{1}^ \top }{\bf{p}} = 1 \} $  and ${\bf{\Omega} }$ is a subset of ${\Delta _N}$.  Agnostic federated learning (AFL)~\citep{mohri2019agnostic} firstly introduces the distributionally robust (agnostic) loss in federated learning and provides the convergence rate for (strongly) convex functions. However, AFL does not discuss the setting of ${\bf{\Omega} }$. DRFA-Prox~\citep{deng2021distributionally} considers ${\bf{\Omega} }={\Delta _N}$ and imposes a regularizer on adversarial distribution to leverage the prior distribution. \textcolor{black}{\textbf{Motivation:} Nonetheless, previous studies have left three pivotal challenges unaddressed, which \textit{motivates} this work to systematically tackle these gaps through a novel optimization framework.} First and foremost is the question of whether it is feasible to establish an uncertainty framework capable of not only flexibly navigating the trade-off between model robustness and performance but also effectively harnessing the information embedded in the prior distribution. Second, addressing the design of asynchronous algorithms with guaranteed convergence is imperative. Compared to synchronous algorithms, the master in asynchronous algorithms can update its parameters after receiving updates from only a small subset of workers~\citep{chang2016asynchronous}.  Asynchronous algorithms are particularly desirable in practice since they can relax strict data dependencies and ensure convergence even in the presence of device failures \citep{zhang2014asynchronous}. Finally, a critical question remains: is it feasible to flexibly adjust the degree of robustness? Furthermore, it is imperative to furnish a convergence guarantee when the objectives (\textit{i.e.}, ${f_j}({\boldsymbol{w}_j}),\forall j$) are non-convex. 

To this end, we introduce ASPIRE-EASE as a means to effectively tackle the aforementioned challenges. Primarily, in contrast to existing methodologies, our formulation integrates the prior distribution into the constraint, which can not only leverage the prior distribution more effectively but also achieve guaranteed feasibility for any adversarial distribution within the uncertainty set. \textcolor{black}{The prior distribution can be obtained based on prior knowledge, which is necessary to construct the uncertainty (ambiguity) set and obtain a more robust model. For instance, when prior knowledge regarding the importance of each node is available, such as through data quality assessments \cite{tang2024adapted}, the prior distribution can be designed to reflect this importance \cite{qian2019robust}. Similarly, if prior knowledge suggests that certain nodes may be malicious, which can be obtained from historical model updates \cite{zhang2022fldetector}, assigning an appropriate prior distribution, e.g., by down-weighting the suspected nodes, can effectively mitigate their impact on the training process \cite{deng2021distributionally}. In the absence of any prior knowledge, a uniform prior distribution can be adopted as a default choice \cite{shapiro2023bayesian}.}  Specifically, we formulate the prior distribution informed distributionally robust optimization (PD-DRO) problem as:
\begin{align}
\label{eq:3}
\mathop {{\rm{min}}}\limits_{\boldsymbol{z}\in{{\boldsymbol{\mathcal{Z}}}},\{{\boldsymbol{w}_j}\in{{\boldsymbol{\mathcal{W}}}}\}} & \mathop {{\rm{max}}} \limits_{{\bf{p}}  \in \boldsymbol{\mathcal{P}}} \sum\nolimits_{j} {{p_j}{f_j}({\boldsymbol{w}_j})} \\ 
{\rm{s.t.}} \; \quad    \boldsymbol{z}&  =    {\boldsymbol{w}_j}, \; j\!=\!1,\!\cdots\!, N , \nonumber\\
{\rm{var.}} \; \;  \; \;  \boldsymbol{z}&,{\boldsymbol{w}_1},{\boldsymbol{w}_2}, \cdots ,{\boldsymbol{w}_N}, \nonumber
\end{align}
where $\boldsymbol{z}\! \in \! {\mathbb{R}^p}$ is the global consensus variable, $\boldsymbol{w}_j\! \in \! {\mathbb{R}^p}$ is the local variable (local model parameter) of $j^{\rm{th}}$ worker and ${\boldsymbol{\mathcal{Z}}}\!\subseteq\! \mathbb{R}^p$ is a nonempty closed convex set.  \textcolor{black}{$\boldsymbol{\mathcal{P}}\!\subseteq\! \mathbb{R}_ + ^N$  is the uncertainty (ambiguity) set of adversarial distribution ${\bf{p}}$. While many existing works adopt the simplex as the ambiguity set, i.e., $\boldsymbol{\mathcal{P}}=\{ {\bf{p}} \in {\mathbb{R}^{N}_{+}}: {\bf{1}^ \top }{\bf{p}} = 1 \}$, incorporating prior distribution can significantly enhance robustness and performance, as discussed in \cite{qian2019robust,deng2021distributionally}. To this end, the ambiguity set $\boldsymbol{\mathcal{P}}$, which is constructed in accordance with the prior distribution, is considered in this work. Specifically, the ambiguity set $\boldsymbol{\mathcal{P}}$ is centered around the prior distribution, with both the component-wise deviations and the total deviation from the prior distribution bounded, which will be discussed in detail in Section \ref{sec:ease}.} {\color{black} Introducing the global consensus variable $\boldsymbol{z}$ in Eq. (\ref{eq:3}) facilitates the development of distributed algorithms for PD-DRO on the parameter-sever architecture \citep{assran2020advances}.}
To solve the PD-DRO problem in an asynchronous distributed manner, we first propose \textbf{A}synchronous \textbf{S}ingle-loo\textbf{P} alternat\textbf{I}ve g\textbf{R}adient proj\textbf{E}ction (ASPIRE), which employs \emph{simple} gradient projection steps for the update of primal and dual variables at every iteration, thus is computationally \emph{efficient}. Next, the it\textbf{E}rative \textbf{A}ctive \textbf{S}\textbf{E}t method (EASE) is employed as a replacement for the traditional cutting plane method. This substitution aims to enhance computational efficiency and expedite convergence. And we furnish a convergence guarantee for the proposed algorithm. Additionally, we introduce an adaptive variant of ASPIRE that allows for flexible adjustment of the number of active workers (\textit{i.e.}, the number of workers that communicate with master at each iteration). Furthermore, a new uncertainty set, \textit{i.e.}, constrained $D$-norm ($CD$-norm), is proposed in this paper and its advantages include: 1) it possesses the inherent capability to flexibly modulate the level of robustness; 2) the resultant subproblem exhibits computational simplicity; 3) it can effectively leverage the prior distribution and flexibly set the bounds for every $p_j$. {\color{black}In addition to the proposed $CD$-norm uncertainty set, we also provide a comprehensive analysis about different uncertainty sets that can be employed in our framework.} The novelty of this work lies in the following aspects: 1) Unlike previous studies, this work incorporates the prior distribution directly as a constraint into the problem formulation, enabling more effective utilization of the prior distribution. 2)  A novel $CD$-norm uncertainty set is introduced, along with a comprehensive analysis comparing it to other commonly used uncertainty sets. 3) In contrast to existing works, this work proposes an asynchronous distributed algorithm with convergence guarantees and an adaptive technique for the FDRO problem.

{\bf{Contributions.}} Our contributions can be summarized as follows:

\begin{enumerate}
    \item We formulate a PD-DRO problem with $CD$-norm uncertainty set. PD-DRO integrates the prior distribution as constraints, enhancing the effective utilization of prior information and ensuring robustness. Additionally, $CD$-norm is developed to delineate the ambiguity set encapsulating the prior distribution, offering a flexible mechanism for regulating the trade-off between model robustness and performance.

    \item We develop a \emph{single-loop} \emph{asynchronous} algorithm, namely ASPIRE-EASE, to optimize PD-DRO in an asynchronous distributed manner. ASPIRE employs simple gradient projection steps to update the variables at every iteration, ensuring computational efficiency. Furthermore, EASE is introduced to replace the cutting plane method, thereby augmenting computational efficiency and expediting convergence. We demonstrate that even if the objectives ${f_j}({\boldsymbol{w}_j}),\forall j$ are non-convex, the proposed algorithm is guaranteed to converge. We also theoretically derive the iteration complexity and communication complexity of ASPIRE-EASE. 

    \item {\color{black}We extend the proposed framework to incorporate a variety of different uncertainty sets, e.g., ellipsoid uncertainty set and Wasserstein-1 distance uncertainty set. We theoretically analyze the computational complexity for each uncertainty set. To accelerate the convergence speed, we further propose ASPIRE-ADP, which can adaptively adjust the number of active workers.}

    \item Extensive empirical investigations conducted on four distinct real-world datasets underscore the superior performance of the proposed algorithm. The results reveal that ASPIRE-EASE not only guarantees the model's robustness in the face of data heterogeneity but also effectively mitigates malicious attacks.

\end{enumerate}

\section{Preliminaries}\label{Preliminaries}

\subsection{Distributionally Robust Optimization}

Optimization problems often contain uncertain parameters. 
A small perturbation of the parameters could render the optimal solution of the original optimization problem infeasible or completely meaningless \citep{bertsimas2004price}. Distributionally robust optimization (DRO) \citep{shi2023distributionally,chen2023adaptive,jia2025distributionally} assumes that the probability distributions of uncertain parameters are unknown but remain in an ambiguity (uncertainty) set and aims to find a decision that minimizes the worst case expected cost over the ambiguity set, whose general form can be expressed as,
\begin{align}
\mathop {\min }\limits_{\boldsymbol{x} \in \boldsymbol{\mathcal{X}}} \mathop {\max }\limits_{P \in {{\bf{P}}}} {\mathbb{E}_P}[r(\boldsymbol{x},\boldsymbol{\xi})],
\end{align}
where $\boldsymbol{x} \! \in \! \boldsymbol{\mathcal{X}} $ represents the decision variable,  ${\bf{P}}$ is the ambiguity set of probability distributions $P$ of uncertain parameters $\boldsymbol{\xi}$. Existing methods for solving DRO can be broadly grouped into two widely-used categories \citep{rahimian2019distributionally}: 1) Dual methods \citep{delage2010distributionally,wiesemann2013robust} reformulate the primal DRO problems as deterministic optimization problems through duality theory. Ben-Tal et al. \citep{ben1999robust} reformulate the robust linear optimization (RLO) problem with an ellipsoidal uncertainty set as a second-order cone optimization problem (SOCP).  2) Cutting plane methods \citep{mehrotra2014cutting,bertsimas2016reformulation} (also called adversarial approaches \citep{gorissen2015practical}) continuously solve an approximate problem with a finite number of constraints of the primal DRO problem, and subsequently check whether new constraints are needed to refine the feasible set. Recently, with the advancement of bilevel optimization, a novel perspective to solving distributionally robust optimization problems has emerged \cite{jiaoasynchronous,hu2023contextual}.


\subsection{Cutting Plane Method for PD-DRO}

In this section, we introduce the cutting plane method for PD-DRO in Eq. (\ref{eq:3}). We first reformulate PD-DRO by introducing an additional variable  $h \!\in\! {{\boldsymbol{\mathcal{H}}}} $ (${{\boldsymbol{\mathcal{H}}}}\! \subseteq\! {\mathbb{R}^1}$ is a nonempty closed convex set) and protection function $g(\{ {\boldsymbol{w}_j}\} )$ \citep{yang2014distributed}. Introducing additional variable $h$  is an epigraph reformulation \citep{yanikouglu2019survey}. In this case, Eq. (\ref{eq:3}) can be reformulated as the form with uncertainty in the constraints:
\begin{align}
 \label{eq:4-1}
\mathop {{\rm{min}}}\limits_{\boldsymbol{z}\in{{\boldsymbol{\mathcal{Z}}}},\{{\boldsymbol{w}_j}\in{{\boldsymbol{\mathcal{W}}}}\}, h \in {{\boldsymbol{\mathcal{H}}}}} &\quad  h  \nonumber\\
{\rm{s.t.}} \; \sum\nolimits_{j}  {\overline{p} {f_j}({\boldsymbol{w}_j})} \!  + \!   g&(\{ {\boldsymbol{w}_j}\} ) - h \! \le \! 0, \\
\boldsymbol{z}  =    {\boldsymbol{w}_j} , \; j\!=&1,\!\cdots\!, N , \nonumber\\
{\rm{var.}}\quad \boldsymbol{z},{\boldsymbol{w}_1},{\boldsymbol{w}_2}, &\cdots ,{\boldsymbol{w}_N}, h. \nonumber
\end{align}
where $\overline{p} $ is the nominal value of the adversarial distribution for every worker and $g(\{ {\boldsymbol{w}_j}\} ) = \mathop {\max }\limits_{{\bf{p}} \in \boldsymbol{\mathcal{P}}} \sum\nolimits_{j} {({p_j} - \overline{p} ){f_j}({\boldsymbol{w}_j})} $ is the protection function. Eq. (\ref{eq:4-1}) is a semi-infinite program (SIP) which contains infinite constraints and cannot be solved directly \citep{rahimian2019distributionally}{}. Denoting the set of cutting plane parameters in $(t\!+\!1)^{\rm{th}}$ iteration as ${{\bf{A}}^t} \!  \subseteq \!  {\mathbb{R}^N}$, the following function is used to approximate $g(\{ {\boldsymbol{w}_j}\} ) $:
\begin{equation}
\begin{aligned}
\overline{g} (\{ {\boldsymbol{w}_j}\} ) = \mathop {\max }\limits_{{\boldsymbol{a}_l} \in {{\bf{A}}^t}} {\boldsymbol{a}^\top_l} {\bf{f}} (\boldsymbol{w}) = \mathop {\max }\limits_{{\boldsymbol{a}_l} \in {{\bf{A}}^t}} \sum\nolimits_{j}  {{a_{l,j}}{{f}_j}({\boldsymbol{w}_j})}, 
\end{aligned}
\end{equation}
where ${\boldsymbol{a}_l = [{a_{l,1}}, \cdots ,{a_{l,N}}]} \! \in \! {\mathbb{R}^N}$ denotes the parameters of $l^{\rm{th}}$ cutting plane in ${{\bf{A}}^t}$ and ${\bf{f}}(\boldsymbol{w}) \! = \! [{{f}_1}(\boldsymbol{w}_1),\cdots, {{f}_N}(\boldsymbol{w}_N)] \! \in \! {\mathbb{R}^N}$. Substituting the protection function $g(\{ {\boldsymbol{w}_j}\} ) $ with $\overline{g} (\{ {\boldsymbol{w}_j}\} )$, we can obtain the following approximate problem:
\begin{align}
\label{eq:9}
\mathop {{\rm{min}}}\limits_{\boldsymbol{z}\in{{\boldsymbol{\mathcal{Z}}}},\{{\boldsymbol{w}_j}\in{{\boldsymbol{\mathcal{W}}}}\}, h \in {{\boldsymbol{\mathcal{H}}}}} & \quad h \nonumber\\
{\rm{s.t.}} \; \sum\nolimits_{j}\! {(\overline{p}  + {a_{l,j}}){{f}_j}({{\boldsymbol{w}}_j})}&  - h \!\le\! 0,{\rm{     }}\forall {\boldsymbol{a}_l} \! \in \! {{\bf{A}}^t}, \\
\boldsymbol{z}  =    {\boldsymbol{w}_j} , \; j\!=&1,\!\cdots\!, N , \nonumber\\
{\rm{var.}}\quad \boldsymbol{z},{\boldsymbol{w}_1},{\boldsymbol{w}_2}, & \cdots ,{\boldsymbol{w}_N}, h. \nonumber
\end{align}

\begin{algorithm}[tb]
   \caption{ASPIRE-EASE}
   \label{algorithm1}
\begin{algorithmic}
 \STATE {\bfseries Initialization:}  iteration $t = 0$, variables $\{{{\boldsymbol{w}}_j^0}\}$, ${\boldsymbol{z}}^0$, ${h}^0$, $\{{\lambda _l^0}\}$, $\{{{\boldsymbol{\phi}}_j^0}\}$  and set ${{\bf{A}}^0} $.
   \REPEAT
   \FOR{\emph{active worker}}
   \STATE updates local  ${\boldsymbol{w}_j^{t+1}}$ according to Eq. (\ref{eq:15});
   \ENDFOR
   \STATE \emph{active workers} transmit local model parameters and loss to \emph{master};
   \STATE \emph{master} receives updates from \emph{active workers}  \textbf{do}
   \STATE  \quad  updates ${\boldsymbol{z}}^{t+1}$, $h^{t+1}$, $\{{\lambda _l^{t + 1}}\}$, $\{{\boldsymbol{\phi}_j^{t+1}}\}$ in master according to Eq. (\ref{eq:16}), (\ref{eq:17}), (\ref{eq:lambda_update}), (\ref{eq:y_update});
    \STATE  \emph{master} broadcasts ${\boldsymbol{z}}^{t+1}$, $h^{t+1}$, $\{{\lambda _l^{t+1}}\}$ to \emph{active workers};
    \FOR{\emph{active worker}}
   \STATE updates local  ${\boldsymbol{\phi}_j^{t+1}}$ according to Eq. (\ref{eq:y_update});
   \ENDFOR
   \IF{$(t+1)$ mod $k$ $==$ 0 and $t< T_1$}
   \STATE   \emph{master} updates ${{\rm \bf{A}}^{t + 1}}$ according to Eq. (\ref{eq:cutting set}) and (\ref{eq:active}), and broadcast parameters to all workers;
   \ENDIF
   \STATE $t =t+1$;
   \UNTIL{convergence} 
\end{algorithmic}
\end{algorithm}

\section{Asynchronous Single-loop Alternative Gradient Projection} 
Distributed optimization emerges as an appealing approach for large-scale learning tasks \citep{yang2008distributed,han2024federated}, as it obviates the need for data aggregation. This not only safeguards data privacy but also diminishes bandwidth requirements \citep{subramanya2021centralized}.  When the neural network models (\textit{i.e.}, ${f_j}({\boldsymbol{w}_j}),\forall j$ are non-convex functions) are used, solving the problem in Eq. (\ref{eq:9}) in a distributed manner faces two challenges. 1) Determining the optimal solution for a non-convex subproblem necessitates a substantial number of iterations, rendering the process highly computationally intensive. Thus, the traditional Alternating Direction Method of Multipliers (ADMM) is ineffective. 2) 
The communication latencies among workers may manifest notable discrepancies \citep{chen2020asynchronous}. As a result, there is a distinct inclination towards the preference for asynchronous algorithms.

To this end, we propose the \textbf{A}synchronous \textbf{S}ingle-loo\textbf{P} alternat\textbf{I}ve g\textbf{R}adient proj\textbf{E}ction (ASPIRE). The advantages of the proposed algorithm include: 1) ASPIRE uses simple gradient projection steps to update variables in each iteration and therefore it is computationally more efficient than the traditional ADMM method, which seeks to find the optimal solution in non-convex (for ${\boldsymbol{w}_j}, \forall j$) and convex (for ${\boldsymbol{z}}$ and $h$) optimization subproblems every iteration,  2) the proposed asynchronous algorithm does not need strict synchronization among different workers. Therefore, ASPIRE remains resilient against communication delays and potential hardware failures from workers. Details of the algorithm are given below. Firstly, we define the node as master which is responsible for updating the global variable $\boldsymbol{z}$, and we define the node which is responsible for updating the local variable ${\boldsymbol{w}_j}$ as worker $j$. In each iteration, the master updates its variables once it receives updates from at least $S$ workers, \textit{i.e.}, active workers, satisfying $1 \le S \le N$. For brevity, we refer to $S$ as the number of active workers (NAW) hereafter. ${{\bf{Q}}^{t + 1}}$ denotes the index subset of workers from which the master receives updates during $(t+1)^{\rm{th}}$ iteration. We also assume the master will receive updated variables from every worker at least once for each $\tau$ iterations. The augmented Lagrangian function of Eq. (\ref{eq:9}) can be written as:
\begin{equation}
\label{eq:11}
\begin{array}{l}
 {L_p}  =  h  +  \sum\nolimits_{l} {{\lambda _l}(\sum\nolimits_{j} {(\overline{p}  + {a_{l,j}}){f_j}({\boldsymbol{w}_j})}  - h)} \\
\qquad   +  \sum\nolimits_{j}\! {{\boldsymbol{\phi}^\top_j}\!(\boldsymbol{z}\! -\! {\boldsymbol{w}_j})} +  \sum\nolimits_{j} \! {\frac{{{\kappa _1}}}{2}||\boldsymbol{z} \! - \! {\boldsymbol{w}_j}|{|^2}}\!,  
\end{array}
\end{equation}
where ${L_p}\!=\!{L_p}{\rm{(\{ }}{\boldsymbol{w}_j}{\rm{\} ,}}\boldsymbol{z},h,\{ {\lambda _l}\} ,\{ {\boldsymbol{\phi}_j}\})$, ${\lambda _l}\!\in\! {{\boldsymbol{\Lambda}}},\forall l$  and ${\boldsymbol{\phi}_j}\! \in \! {{\boldsymbol{\Phi}}},\forall j$ represent the dual variables of inequality and equality constraints in Eq. (\ref{eq:9}), respectively. ${{\boldsymbol{\Lambda}}}\! \subseteq \! \mathbb{R}^1$ and ${{\boldsymbol{\Phi}}}\! \subseteq \! \mathbb{R}^p$ are nonempty closed convex sets, constant ${\kappa _1} > 0$  is a penalty parameter. Note that Eq. (\ref{eq:11}) does not consider the second-order penalty term for inequality constraint since it will invalidate the distributed optimization. Following \citep{xu2020unified}, the regularized version of Eq. (\ref{eq:11}) is employed to update all variables as follows,
\begin{align}
\! {\widetilde{L}_p}{\rm{(\{ }}{\boldsymbol{w}_j}{\rm{\} ,}}\boldsymbol{z},\!h,\!\{ {\lambda _l}\} ,\!\{ {\boldsymbol{\phi}_j}\}) 
 \!=\!{L_p} \!-\! \sum\nolimits_{l} \!{\frac{{c_1^{t}}}{2}||{\lambda _l}|{|^2}}  \!-\! \sum\nolimits_{j}\! {\frac{{c_2^{t}}}{2}||{\boldsymbol{\phi}_j}|{|^2}},
\end{align}
where $c_1^{t}$ and $c_2^{t}$ denote the regularization terms in $(t+1)^{\rm{th}}$ iteration.  In $(t+1)^{\rm{th}}$ master iteration, the proposed algorithm proceeds as follows.

\textbf{1)} \emph{Active} \emph{workers} update the local variables ${\boldsymbol{w}_j}$ as follows,
\begin{equation}
\label{eq:15}
{\boldsymbol{w}_j^{t+1}} = \left\{ \begin{array}{l}
{\mathcal{P}_{{\boldsymbol{\mathcal{W}}}}}({\boldsymbol{w}_j^{t}}  -  {\alpha _{\boldsymbol{w}}^{\widetilde{t}_j}}{\nabla _{{\boldsymbol{w}_j}}}{\widetilde{L}_p^{\widetilde{t}_j}}),\forall j  \in  {{\bf{Q}}^{t + 1}},\\
{\boldsymbol{w}_j^t},\forall j \notin {{\bf{Q}}^{t + 1}},
\end{array} \right.
\end{equation}
where  $\widetilde{t}_j$ is the last iteration during which local worker $j$ was active, and ${\widetilde{L}_p^{\widetilde{t}_j}}$ is the simplified form of $ {\widetilde{L}_p}{\rm{(\{ }}{\boldsymbol{w}_j^{\widetilde{t}_j}}{\rm{\} ,}}\boldsymbol{z}^{\widetilde{t}_j},h^{\widetilde{t}_j},\!\{ {\lambda _l^{\widetilde{t}_j}}\} ,\!\{ {\boldsymbol{\phi}_j^{\widetilde{t}_j}}\} {\rm{)}}$, we have that,
{\color{black}
\begin{equation}
\label{eq:new_12}
\begin{array}{l}
    {\nabla _{{\boldsymbol{w}_j}}}{ \widetilde{L}_p}{\rm{(\{ }}{\boldsymbol{w}_j^{\widetilde{t}_j}}{\rm{\} ,}}\boldsymbol{z}^{\widetilde{t}_j},h^{\widetilde{t}_j},\!\{ {\lambda _l^{\widetilde{t}_j}}\} ,\!\{ {\boldsymbol{\phi}_j^{\widetilde{t}_j}}\} {\rm{)}} 
    \\
    = \sum\nolimits_{l}\! {{\lambda _l^{\widetilde{t}_j}}  {(\overline{p}  + {a_{l,j}}){\nabla _{{\boldsymbol{w}_j}}}{f_j}({\boldsymbol{w}_j^{\widetilde{t}_j}})}} - {\boldsymbol{w}_j^{\widetilde{t}_j}}  +    {\kappa _1}({\boldsymbol{w}_j^{\widetilde{t}_j}}-\boldsymbol{z}^{\widetilde{t}_j}). 
\end{array}
\end{equation}

It is seen from Eq. (\ref{eq:15}) and (\ref{eq:new_12}) that the local update of $\boldsymbol{w}_j$ is only based on the variables on the local worker $j$ and master.} And $ {\boldsymbol{w}_j^{t}} = {\boldsymbol{w}_j^{\widetilde{t}_j}} $, ${\boldsymbol{\phi}_j^{t}} = {\boldsymbol{\phi}_j^{\widetilde{t}_j}},\forall j \! \in \! {{\bf{Q}}^{t + 1}}$, ${\alpha _{\boldsymbol{w}}^{\widetilde{t}_j}}$ represents the step-size and we set  ${\alpha _{\boldsymbol{w}}^t} = {\eta _{\boldsymbol{w}}^t}$ when $t<T_1$ and ${\alpha _{\boldsymbol{w}}^t} = {\underline{\eta _{\boldsymbol{w}}}}$ when $t \ge T_1$, where  ${\eta _{\boldsymbol{w}}^t}$ and constant ${\underline{\eta _{\boldsymbol{w}}}}$ will be introduced below. $\mathcal{P}_{{\boldsymbol{\mathcal{W}}}}$ represents the projection onto the closed convex set ${{\boldsymbol{\mathcal{W}}}}$ and we set ${{\boldsymbol{\mathcal{W}}}} = \{ {{{\boldsymbol{w}}_j}}|\;|| {{\boldsymbol{w}}_j}||_{\infty} \! \le {\alpha _1}\} $, ${\alpha _1}$ is a positive constant. And  
then, the active workers ($ j \! \in \! {{\bf{Q}}^{t + 1}}$) transmit their local model parameters ${{\boldsymbol{w}}_j^{t+1}}$ and loss ${f}_j({\boldsymbol{w}_j})$ to the master.

\textbf{2)} After receiving the updates from active workers, the \emph{master} updates the global consensus variable $\boldsymbol{z}$, additional variable $h$ and dual variables ${\lambda _l}$ as follows,
\begin{equation}
\label{eq:16}
\boldsymbol{z}^{t+1}\! =\! {\mathcal{P}_{{\boldsymbol{\mathcal{Z}}}}}({\boldsymbol{z}}^{t} - {\eta _{\boldsymbol{z}}^t}{\nabla _{\boldsymbol{z}}}{\widetilde{L}_p}{\rm{(\{ }}{{\boldsymbol{w}}_j^{t+1}}{\rm{\} ,}}{\boldsymbol{z}}^{t},h^{t},\! \{ {\lambda _l^{t}}\} ,\! \{ {{\boldsymbol{\phi}}_j^{t}}\} {\rm{)}}),
\end{equation}
\begin{equation}
\label{eq:17}
h^{t+1} \!=\! {\mathcal{P}_{{\boldsymbol{\mathcal{H}}}}}(h^{t} - {\eta _h^t}{\nabla _h}{\widetilde{L} _p}{\rm{(\{ }}{{\boldsymbol{w}}_j^{t+1}}{\rm{\} ,}}{\boldsymbol{z}}^{t+1},h^{t},\! \{ {\lambda _l^{t}}\} ,\!\{ {{\boldsymbol{\phi}}_j^{t}}\} {\rm{)}}),
\end{equation}
\begin{equation}
\label{eq:lambda_update}
{\lambda _l^{t+1}} \!=\! {\mathcal{P}_{{\boldsymbol{\Lambda}}} }({\lambda _l^{t}} \! + \! {\rho _1}{\nabla _{{\lambda _l}}}{\widetilde{L}_p}{\rm{(\{ }}{{\boldsymbol{w}}_j^{t+1}}{\rm{\} ,}}{\boldsymbol{z}}^{t+1},h^{t+1},\! \{ {\lambda _l^{t}}\} , \! \{ {{\boldsymbol{\phi}}_j^{t}}\} )),
\end{equation}
where ${\eta _{\boldsymbol{z}}^t}$, ${\eta _h^t}$ and $\rho_1$ represent the step-sizes. ${\mathcal{P}_{{\boldsymbol{\mathcal{Z}}}}}$, ${\mathcal{P}_{{\boldsymbol{\mathcal{H}}}}}$ and ${\mathcal{P}_{{\boldsymbol{\Lambda}}}}$ respectively represent the projection onto the closed convex sets ${{\boldsymbol{\mathcal{Z}}}}$, ${{\boldsymbol{\mathcal{H}}}}$ and ${{\boldsymbol{\Lambda}}}$. We set ${{\boldsymbol{\mathcal{Z}}}} = \{ {{{\boldsymbol{z}}}}|\;|| {{\boldsymbol{z}}}||_{\infty} \! \le {\alpha _1}\} $, ${{\boldsymbol{\mathcal{H}}}} = \{ h|\; 0 \le \! h \! \le {\alpha _2}\} $ and ${{\boldsymbol{\Lambda}}} = \{ {\lambda _l}|\; 0 \le \! {\lambda _l} \! \le {\alpha _3}\} $, where $\alpha _2$ and $\alpha _3$ are positive constants. $|{{\bf{A}}^t}|$ denotes the number of cutting planes. Then, master broadcasts ${\boldsymbol{z}}^{t+1}$, $h^{t+1}$, $\{{\lambda _l^{t+1}}\}$ to the active workers.

\textbf{3)} \emph{Active} \emph{workers} update the local dual variables ${\boldsymbol{\phi}_j}$ as follows,
\begin{equation}
\label{eq:y_update}
{\boldsymbol{\phi}_j^{t+1}}\! =\! \left\{\! \begin{array}{l}
{\mathcal{P}_{{\boldsymbol{\Phi}}}}({{\boldsymbol{\phi}}_j^{t}}\! +\! {\rho _2}{\nabla _{{{\boldsymbol{\phi}}_j}}}{\nabla _{{{\boldsymbol{\phi}}_j}}}g^t{\rm{)}},\forall j \!\in\! {{\bf{Q}}^{t + 1}},\\
{{\boldsymbol{\phi}}_j^{t}},\forall j \notin {{\bf{Q}}^{t + 1}},
\end{array} \right.
\end{equation}
where ${\nabla _{{{\boldsymbol{\phi}}_j}}}g^t={\nabla _{{{\boldsymbol{\phi}}_j}}}{\widetilde{L}_p}{\rm{(\{ }}{{\boldsymbol{w}}_j^{t+1}}{\rm{\} ,}}{\boldsymbol{z}}^{t+1},h^{t+1},\! \{ {\lambda _l^{t+1}}\} ,\! \{ {{\boldsymbol{\phi}}_j^{t}}\} {\rm{)}}$, $\rho_2$ represents the step-size and  ${\mathcal{P}_{{\boldsymbol{\Phi}}}}$ represents the projection onto the closed convex set  ${{\boldsymbol{\Phi}}}$ and we set ${{\boldsymbol{\Phi}}} = \{ {{{\boldsymbol{\phi}}_j}}|\;|| {{\boldsymbol{\phi}}_j}||_{\infty} \! \le {\alpha _4}\} $, $\alpha _4$ is a positive constant. And master can also obtain $\{{\boldsymbol{\phi}_j^{t+1}}\}$ according to Eq. (\ref{eq:y_update}). It is seen that the projection operation in each step is computationally simple since the closed convex sets have simple structures \citep{bertsekas1997nonlinear}.

\section{Iterative Active Set Method}\label{sec:ease}

Cutting plane methods may give rise to numerous linear constraints and lots of extra message passing \citep{yang2014distributed}. To enhance computational efficiency and expedite convergence, we contemplate the removal of inactive cutting planes. The proposed it\textbf{E}rative \textbf{A}ctive \textbf{S}\textbf{E}t method (EASE) can be divided into the two steps: during $T_1$ iterations, 1) solving the cutting plane generation subproblem to generate cutting plane, and 2) removing the inactive cutting plane every $k$ iterations, where $k\!>\!0$ is a pre-set constant and can be controlled flexibly.

The cutting planes are generated based on the specific uncertainty set. For instance, in the case of employing an ellipsoid uncertainty set, the cutting plane is generated by solving a SOCP. In this paper, we propose  $CD$-norm uncertainty set, which is formulated as follows:
\begin{align}
\label{eq:D norm}
\boldsymbol{\mathcal{P}} \!= \! \{ {\bf{p}} \!:  - \widetilde{p}_j \! \le \!  p_j -q_j \! \le \! \widetilde{p}_j,   \sum\nolimits_{j} \!|{\frac{{  p_j -q_j }}{ \widetilde{p}_j  }}| \! \le \! \Gamma   ,  {\bf{1}^ \top }{\bf{p}}\! =\! 1  \},
\end{align}
where  $\Gamma \! \in \! {\mathbb{R}^1}$ can flexibly control the level of robustness, ${\bf{q}}=[q_1, \cdots, q_N] \! \in \! {\mathbb{R}^N}$ represents the {\color{black}\emph{prior} \emph{distribution}}, $-\widetilde{p}_j$ and $\widetilde{p}_j$ ($\widetilde{p}_j \ge 0$) represent the lower and upper bounds for $p_j -q_j$, respectively. The setting of ${\bf{q}}$ and $\widetilde{p}_j, \forall j$ are based on the prior knowledge. $D$-norm is a classical uncertainty set (which is also called as budget uncertainty set) \citep{bertsimas2004price}. We call Eq. (\ref{eq:D norm}) $CD$-norm uncertainty set since ${\bf{p}}$ is a probability vector so all the entries of this vector are non-negative and add up to exactly one, \textit{i.e.}, ${\bf{1}^ \top }{\bf{p}} = 1$. Due to the special structure of $CD$-norm, the cutting plane generation subproblem is easy to solve and the level of robustness in terms of the outage probability, \textit{i.e.}, probabilistic bounds of the violations of constraints can be flexibly adjusted via a single parameter $\Gamma$. We claim that $l_1$-norm (or twice total variation distance) uncertainty set is closely related to $CD$-norm uncertainty set. Nevertheless, there are two differences: 1) $CD$-norm uncertainty set could be regarded as a weighted $l_1$-norm with additional constraints. 2) $CD$-norm uncertainty set can flexibly set the lower and upper bounds for every $p_j$ (\textit{i.e.}, $q_j \! - \! \widetilde{p}_j  \! \le \!  p_j \! \le \! p_j \! + \! \widetilde{p}_j$), while $0 \! \le \! p_j  \! \le \! 1, \forall j$ in $l_1$-norm uncertainty set. Based on the $CD$-norm uncertainty set, the cutting plane can be derived as follows,

1) Solve the following problem,
\begin{align}
\label{eq:18}
& {\bf{p}}^{t + 1}  = \mathop {\arg \max }\limits_{  {p_1},\cdots,{p_N} }  \sum\nolimits_{j} {({p_j}  - \overline{p}  ){{ f}_j}({{\boldsymbol{w}_j}})} \nonumber \vspace{-3mm}\\
{\rm{s.t.}}\;  \sum\nolimits_{j} & |{\frac{{  p_j \!- \!q_j }}{ \widetilde{p}_j  }}|  \! \le \!  \Gamma,\;  - \widetilde{p}_j \! \le \!   p_j \!- \! q_j  \! \le \!  \widetilde{p}_j, \forall{j},\; \sum\nolimits_{j}\! {p_j} \! = \! 1  \\
&{\rm{var.}}\quad \quad \quad \quad \; {p_1},\cdots,{p_N},  \nonumber
\end{align}
where ${\bf{p}}^{t + 1} \! = \! [ p_1^{t + 1},\! \cdots, p_N^{t + 1}]\! \in \! {\mathbb{R}^N} $. Let $\widetilde{{\bf{a}}}^{t + 1} \!= \! {\bf{p}} ^{t + 1} - \overline{{\bf{p}}}  $, where $\overline{{\bf{p}}}    = [\overline{p}  , \cdots ,\overline{p} ] \! \in \! {\mathbb{R}^N}$. This first step aims to obtain the distribution $\widetilde{{\bf{a}}}^{t+1}$  by solving problem in Eq. (\ref{eq:18}). This problem can be effectively solved through combining merge sort \citep{cole1988parallel,jiao2025pr} (for sorting $\widetilde{p}_j{f}_j({{\boldsymbol{w}_j}}), j\!=\!1,\cdots,N$) with few basic arithmetic operations (for obtaining $p_j^{t + 1}, j\!=\!1,\cdots,N$). Since $N$ is relatively large in distributed system, the arithmetic complexity of solving problem in Eq. (\ref{eq:18}) is dominated by merge sort, which can be regarded as $\mathcal{O}(N\log (N))$.
 
2) Let ${\bf{f}}(\boldsymbol{w}) \! =\! [{{f}_1}(\boldsymbol{w}_1),\cdots , {{f}_N}(\boldsymbol{w}_N)]\!\in\!{\mathbb{R}^N}$, check the feasibility of the following constraints:
\begin{equation}
\label{eq:19}
{\widetilde{{\bf{a}}}^{t+1}}{}^\top {\bf{f}}({\boldsymbol{w}}) \! \le \!\mathop {\max }\limits_{{\boldsymbol{a}_l} \in {{\bf{A}}^t}} {\boldsymbol{a}_l}{}^\top {\bf{f}} (\boldsymbol{w}),
\end{equation}

\renewcommand\arraystretch{1.2}
\renewcommand\tabcolsep{5pt}
\begin{table}[t]
\centering
\caption{\textcolor{black}{Comparisons between this work and related work.}}
\begin{tabular}{l|c|c|c|c|c|c}
\toprule
  & \cite{qian2019robust} & \cite{mohri2019agnostic} & \cite{issaid2021federated} &  \cite{shen2024stochastic} & \cite{deng2021distributionally} & This work \\ \hline

Non-convex & \checkmark &  & \checkmark &  \checkmark &  \checkmark &  \checkmark \\

Federated &  & \checkmark & \checkmark & \checkmark &  \checkmark &  \checkmark \\

Prior Distribution &  & &  & &  \checkmark &  \checkmark \\

Flexibility &  & & &  & &  \checkmark \\

Asynchronous &  & &  &  & &  \checkmark \\
\bottomrule
\end{tabular}
\label{tab:research_gap}
\end{table}

3) If Eq. (\ref{eq:19}) is violated,  $\widetilde{{\bf{a}}}^{t+1}$ will be added into ${{\rm \bf{A}}^t}$:
\begin{equation}
\label{eq:cutting set}
{{\rm \bf{A}}^{t + 1}} = \left\{ \begin{array}{l}
{{\rm \bf{A}}^t} \cup \{ \widetilde{{\bf{a}}}^{t+1} \},{\rm{ if \;  Eq. (\ref{eq:19}) \; is \; violated }},\\
{{\rm \bf{A}}^t},{\rm{ otherwise}},
\end{array} \right.
\end{equation}
when a new cutting plane is added, its corresponding dual variable ${\lambda _{|{{\bf{A}}^t}| + 1}^{t+1}}=0$ will be generated. After the cutting plane subproblem is solved, the inactive cutting plane will be removed, that is:
\begin{equation}
\label{eq:active}
{{\rm \bf{A}}^{t + 1}} =\left\{ \begin{array}{l}
\complement_{ {{\rm \bf{A}}^{t+1}}} \{{\boldsymbol{a}_l}\} ,{\rm{if}} \; \lambda_l^{t+1}\!=\!0 \, {\rm{and}}\, \lambda_l^{t} \!=\! 0, 1\!\le\!l \!\le\! |{{\bf{A}}^t}|, \\
{{\rm \bf{A}}^{t+1}},{\rm{ otherwise}},
\end{array} \right.
\end{equation}
where $\complement_{ {{\rm \bf{A}}^{t+1}}}\{{\boldsymbol{a}_l}\}$ is the complement of $\{{\boldsymbol{a}_l}\}$ in ${ {{\rm \bf{A}}^{t+1}}}$, and the dual variable will be removed. Then master broadcasts ${{\rm \bf{A}}^{t + 1}}$, $\{ \lambda_l^{t+1} \}$ to all workers.  Details of the proposed algorithm are summarized in Algorithm \ref{algorithm1}. In addition, we provide a comparison between the proposed method and existing approaches across several key aspects, including the non-convexity of the objective function, applicability to federated learning settings, the ability to incorporate prior distribution, flexibility in adjusting the degree of robustness, and the capacity to handle the asynchronous updates. The comparison results are summarized in Table \ref{tab:research_gap}.

\begin{table*}[t]
\caption{Comparison among different uncertainty sets.}
\renewcommand\arraystretch{2.2}
\renewcommand\tabcolsep{3pt}
\label{tab:uncertainty set comparison}
\centering
\scalebox{0.95}{\begin{tabular}{l|c|c|c|c} 
\toprule
Uncertainty Set &  $\boldsymbol{\mathcal{P}}$  & Formulation  & Flexibility$^2$  & Complexity$^3$
\\
\hline
Box & $\{ {\bf{p}}:{p_j^{low}} \le {p_j} \le {\rm{ }}{p_j^{upp}},{\bf{1}^ \top }{\bf{p}} = 1\}$  & $\rm{\underline{LP}}^1$   &  2$N$  & $\mathcal{O}(n\log (n))$ \\
\hline
Ellipsoid &  $\{ {\bf{p}}:{({\bf{p}} - {\bf{q}})^T}{{\bf{Q}}^{ - 1}}({\bf{p}} - {\bf{q}}) \le {\beta},{\bf{1}^ \top }{\bf{p}} = 1\}$ &  SOCP  & 1 & $\mathcal{O}((m\!+\!1)\!^{1\!/\!2}n({n^{2}\!+\!m\!+\!\sum\limits_{i = 1}^{m} \!{k_i^2} })\! \log (\frac{1}{{\varepsilon '}}))$  \\
\hline
Polyhedron & $ \{ {\bf{p}}: {\bf{D}}{\bf{p}} \preceq {\bf{c}},{\bf{1}^ \top }{\bf{p}} = 1\}$   &  $\rm{{LP}}$  & ${L_{\rm{in}} \! \times \! ({N\!+\!1})}$ & $\mathcal{O}((m\!+\!n)^{3/2}{n^{2}}\log (\frac{1}{{\varepsilon '}}))$  \\
\hline
KL-Divergence &  $ \{ {\bf{p}}: \sum\limits_{j = 1}^N {{p_j}\log \frac{{{p_j}}}{{{q_j}}}}   \le {\beta},{\bf{1}^ \top }{\bf{p}} = 1\}$  &  REP  &  1 & $\mathcal{O}({(n)^{7/2}}|\log (\varepsilon ')|)$ \\
\hline
Wasserstein-1 Distance & $\{ {\bf{p}}\!:\! \mathop {\min }\limits_{\gamma  \in \Pi ({\bf{p}},{\bf{q}})} \sum\limits_{i = 1}^N \! {\sum\limits_{j = 1}^N \! {\gamma ({x_i},{y_j})||{x_i} - {y_j}||} } \! \le \! {\beta},{\bf{1}^ \top }{\bf{p}}\! =\! 1\}$   &  $\rm{LP}$  & 1  & $\mathcal{O}((m\!+\!n)^{3/2}{n^{2}}\log (\frac{1}{{\varepsilon '}}))$ \\
\hline
$CD$-norm   & $\{ {\bf{p}} \!: \! - \widetilde{p}_j \! \le \!  p_j -q_j \! \le \! \widetilde{p}_j, \!  \sum\limits_{j = 1}^N \!|{\frac{{  p_j -q_j }}{ \widetilde{p}_j  }}| \! \le \! \Gamma   ,  {\bf{1}^ \top }{\bf{p}} = 1  \}$   &  $\rm{\underline{LP}}$    & 1 & $\mathcal{O}(n\log (n))$ \\
\bottomrule
\end{tabular}}
\\\footnotesize{{\color{black}$^1$ $\rm{\underline{LP}}$ represents the Linear Programming which can be solved by combining merge sort with a few basic arithmetic operations.  $^2$ Flexibility denotes the number of parameters that are utilized to tradeoff between robustness and performance, a lower value indicates better flexibility. $^3$ Complexity denotes the arithmetic complexity of solving the cutting plane generation subproblem.}}
\end{table*}

\begin{figure}[t]
\subfigure[$m=1$, $\varepsilon'=0.9$]
{\begin{minipage}{4.2cm}
    \includegraphics[scale=0.27]{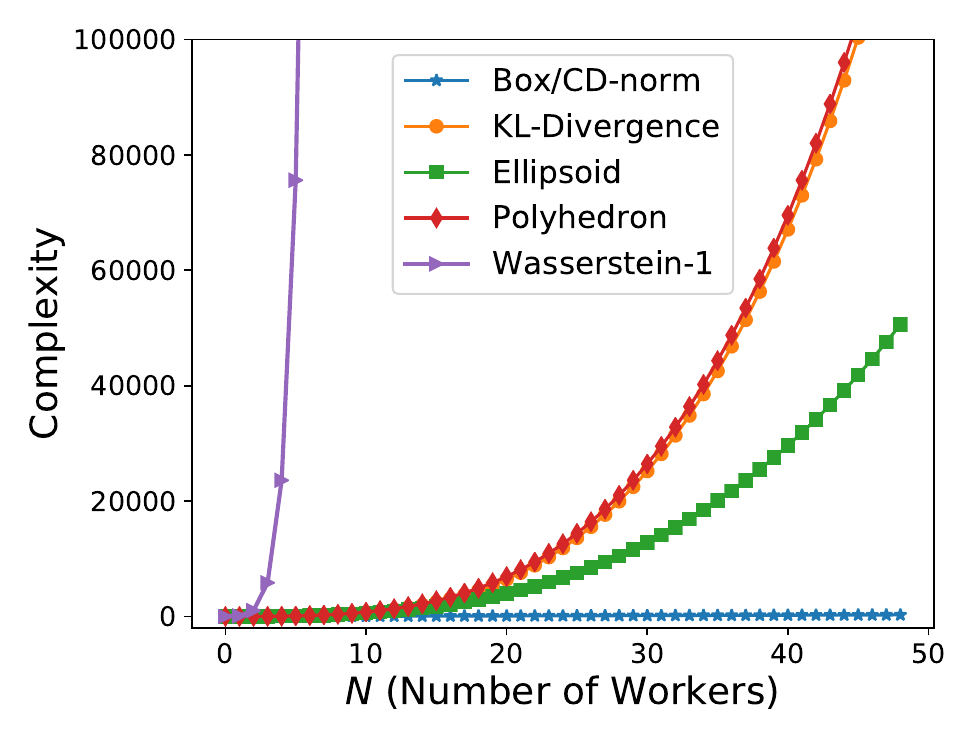}  
\end{minipage}}
\subfigure[$m=5$, $\varepsilon'=0.9$] 
{\begin{minipage}{4.2cm}
    \includegraphics[scale=0.27]{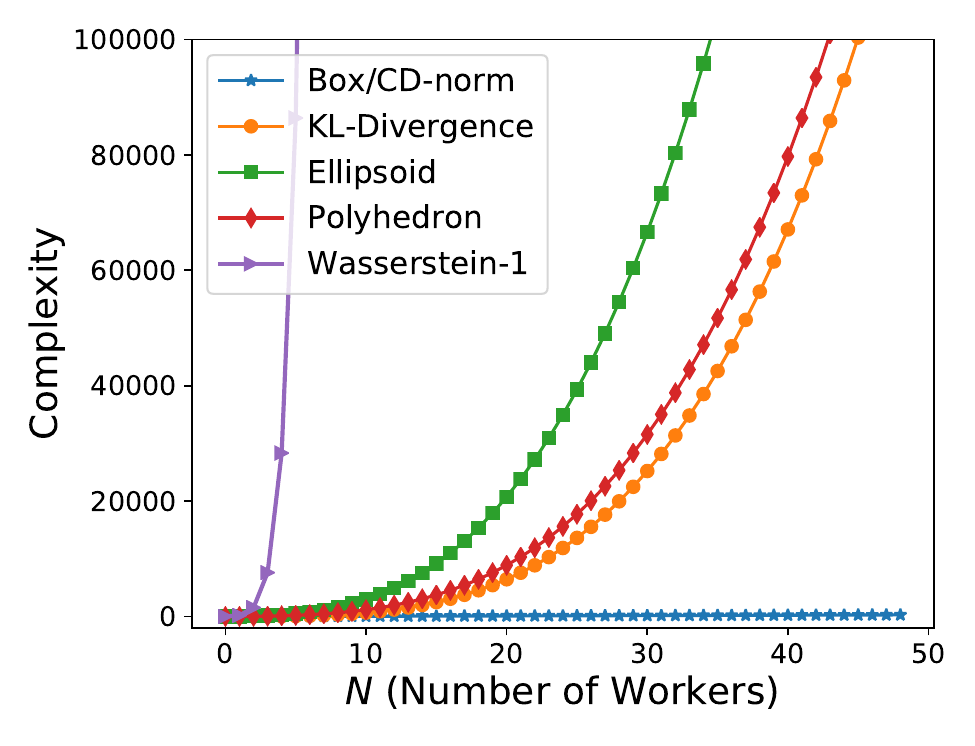}  
\end{minipage}}

\caption{{\color{black}The arithmetic complexity of solving the cutting plane generation subproblem when utilizing different uncertainty sets when setting $\varepsilon'=0.9$. Different numbers of constraints are considered in the subfigures: in subfigure (a), we set $m=1$, while in subfigure (b), we set $m=5$.}} 
\label{fig:complexity}
\vspace{-2mm}
\end{figure}

\begin{figure}[t]
\subfigure[$m=1$, $\varepsilon'=0.1$]
{\begin{minipage}{4.2cm}
    \includegraphics[scale=0.27]{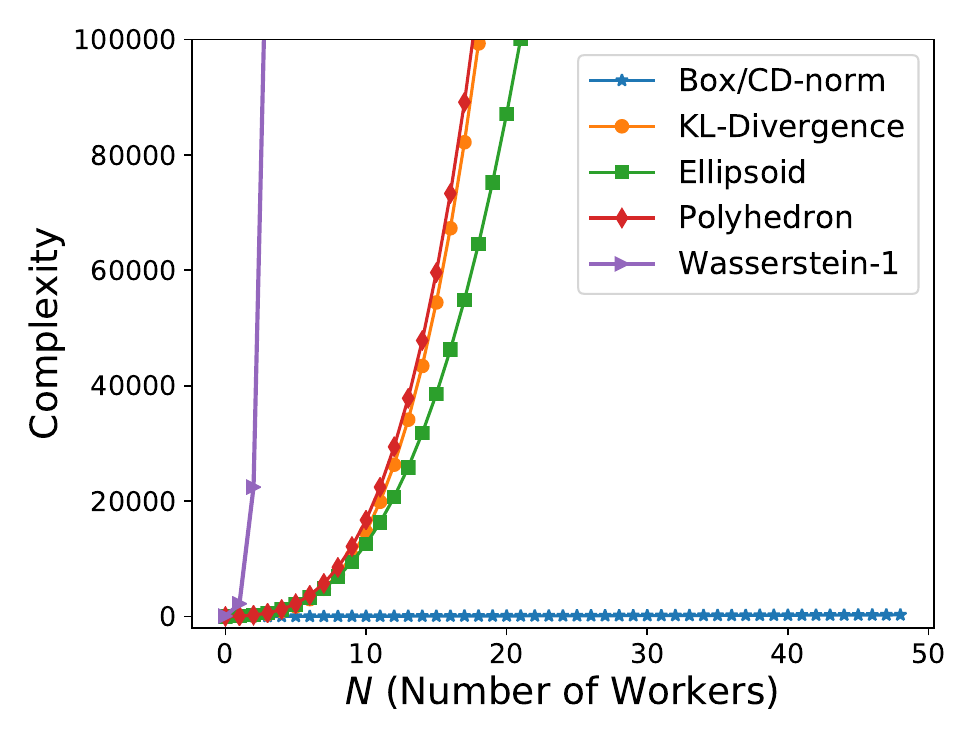}  
\end{minipage}}
\subfigure[$m=5$, $\varepsilon'=0.1$] 
{\begin{minipage}{4.2cm}
    \includegraphics[scale=0.27]{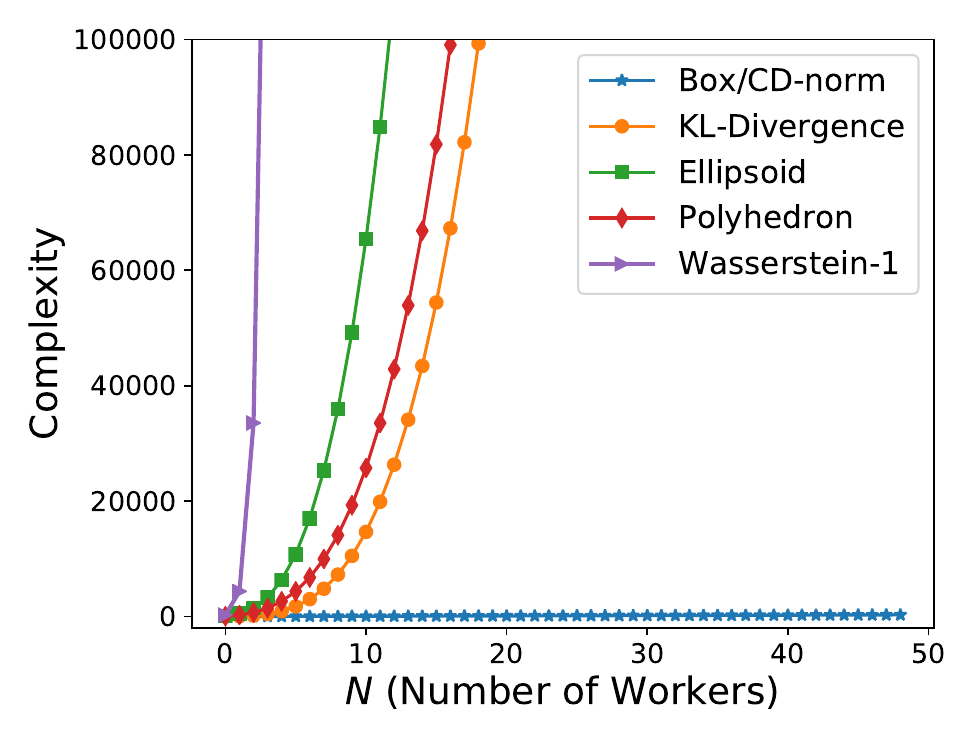}  
\end{minipage}}

\caption{{\color{black}The arithmetic complexity of solving the cutting plane generation subproblem when utilizing different uncertainty sets when setting a relatively smaller $\varepsilon'$, i.e., $\varepsilon'=0.1$. Different numbers of constraints are considered in the subfigures: in subfigure (a), we set $m=1$, while in subfigure (b), we set $m=5$.}} 
\label{fig:complexity}
\vspace{-3mm}
\end{figure}

{\color{black}
\section{Uncertainty Sets}
\label{Uncertainty Sets}

\textcolor{black}{The proposed $CD$-norm uncertainty set offers notable advantages, such as flexibility and computational efficiency. Beyond the $CD$-norm, a variety of alternative uncertainty sets can also be incorporated into our framework, each providing distinct benefits. For instance, ellipsoid uncertainty sets are typically less conservative, leading to robust solutions that do not significantly compromise optimality \cite{bertsimas2004price}, while still providing sufficient coverage of plausible perturbations to balance robustness and performance effectively. Uncertainty sets based on the  KL-divergence and the Wasserstein distance often provide enhanced statistical robustness \cite{namkoong2016stochastic,liu2021stable}.} In this section, different uncertainty sets that can be used in our framework are discussed. Specifically, we respectively formulate the cutting plane generation subproblems for different uncertainty sets, and analyze the arithmetic complexity of solving these subproblems. Moreover, we also focus on whether utilizing different uncertainty sets can flexibly tradeoff between robustness with
performance. \textcolor{black}{The discussions on how to select different uncertainty sets based on various scenarios are provided in Sec. F.}

\subsection{Box Uncertainty Set}
Box uncertainty set was proposed in \citep{soyster1973convex}, which utilizes the box to characterize the uncertainty set and can be written as,
\begin{align}
\label{eq:Box}
\boldsymbol{\mathcal{P}} =  \{ {\bf{p}}:{p_j^{low}} \le {p_j} \le {\rm{ }}{p_j^{upp}},{\bf{1}^ \top }{\bf{p}} = 1\}.
\end{align}
where ${p_j^{low}}$ and ${p_j^{upp}}, \forall j$ are preset constants. The interval of every uncertain coefficient is specified by the box uncertainty set, \textit{i.e.}, ${p_j^{low}} \le {p_j} \le {\rm{ }}{p_j^{upp}}$. When utilizing the box uncertainty set, the following cutting plane generation subproblem is required to be solved in the process of updating cutting planes,
\begin{align}
\label{eq:Box_cutting}
& {\bf{p}}^{t+1}  = \mathop {\arg \max }\limits_{  {p_1},\dots,{p_N} }  \sum\nolimits_{j = 1}^N {({p_j}  - \mathop {{\rm{ }}p}\limits^ -  ){{f}_j}({{\boldsymbol{w}_j}})}\\
&{\rm{s.t.}}\quad \quad  {p_j^{low}} \le {p_j} \le {\rm{ }}{p_j^{upp}} ,\sum\nolimits_{j = 1}^N {{p_j} }  = 1,  \nonumber\\
&{\rm{var.}}\quad \quad \quad \quad \; {p_1},\dots,{p_N}.  \nonumber
\end{align}

It is seen that the problem in Eq. (\ref{eq:Box_cutting}) is an LP. Similar to the problem in Eq. (\ref{eq:18}), Eq. (\ref{eq:Box_cutting}) can be efficiently solved through combining merge sort (for sorting ${f}_j({{\boldsymbol{w}_j}}), j\!=\!1,\dots,N$) with few basic arithmetic operations (for obtaining $p_j(t\!+\!1), j\!=\!1,\dots,N$). As mentioned before, the arithmetic complexity is  $\mathcal{O}(n\log (n))$, where $n=N$ in this problem. Nevertheless, the box uncertainty set is generally considered to be too conservative, which tends to induce over-conservative decisions \citep{bertsimas2004price,li2011comparative}. And the box uncertainty set cannot flexibly tradeoff between robustness with performance since it is required to adjust 2$N$ parameters (\textit{i.e.}, lower and upper bounds for every $p_j$).

\subsection{Ellipsoid Uncertainty Set}
We next proceed with discussions on the  ellipsoid uncertainty set. Firstly, the ellipsoid uncertainty set is given by,
\begin{align}
\label{eq:ellipsoid}
\boldsymbol{\mathcal{P}} = \{ {\bf{p}}:{({\bf{p}} - {\bf{q}})^T}{{\bf{Q}}^{ - 1}}({\bf{p}} - {\bf{q}}) \le {\beta},{\bf{1}^ \top }{\bf{p}} = 1\},
\end{align}
where ${\bf{Q}}\!\in\!{\mathbb{R}^{N \times N}}$ is a positive definite matrix. And the ellipsoid is a ball when ${\bf{Q}}$ is an identity matrix (\textit{i.e.}, ${\bf{Q}}={\bf{I}}\in{\mathbb{R}^{N \times N}}$).
The ellipsoid is widely employed to approximate complicated uncertainty sets since it can succinctly describe a set of discrete points in Euclidean geometry \citep{yang2014distributed,boyd2004convex}. Compared with the box uncertainty set, the ellipsoid uncertainty set is less conservative, but more computationally more intensive \citep{bertsimas2004price} since it leads to a nonlinear optimization subproblem, as given below.
\begin{align}
\label{eq:ellipsoid_cutting}
& {\bf{p}}^{t+1}  = \mathop {\arg \max }\limits_{{p_1},\dots,{p_N} }  \sum\nolimits_{j = 1}^N {({p_j}  - \mathop {{\rm{ }}p}\limits^ -  ){{f}_j}({{\boldsymbol{w}_j}})}\\
&{\rm{s.t.}}\quad  {({\bf{p}} - {\bf{q}})^T}{{\bf{Q}}^{ - 1}}({\bf{p}} - {\bf{q}}) \le {\beta},\sum\nolimits_{j = 1}^N {{p_j} }  = 1,  \nonumber\\
&{\rm{var.}}\quad \quad \quad \quad \; {p_1},\dots,{p_N}.  \nonumber
\end{align}

It is seen that the problem in Eq. (\ref{eq:ellipsoid_cutting}) is a SOCP, which can be solved in polynomial time through the interior point method \citep{boyd2004convex}. Specifically, the arithmetic complexity of interior point method to find the $\varepsilon'$-solution for SOCP is upper bounded by $\mathcal{O}((m\!+\!1)^{1/2}n({n^{2}\!+\!m\!+\!\sum\limits_{i = 1}^{m} {k_i^2} })\log (\frac{1}{{\varepsilon '}}))$ \citep{ben2011lectures}, where $m$ and $n$ are respectively the number of inequality constraints and variables, and $k_i$ represents that the $i$-th inequality constraint is a $k_i+1$ dimension second-order cone. In this problem, $n=N$, $m=1$, $k_1=N$. Compared with the box uncertainty set, the ellipsoid uncertainty set can flexibly tradeoff between robustness with performance by adjusting a single parameter $\beta$.

\subsection{Polyhedron Uncertainty Set}
The form of polyhedron uncertainty set is given by,
\begin{align}
\label{eq:polyhedron}
\boldsymbol{\mathcal{P}} =  \{ {\bf{p}}: {\bf{D}}{\bf{p}} \preceq {\bf{c}}, {\bf{1}^ \top }{\bf{p}} = 1\},
\end{align}
where ${\bf{D}} \in {\mathbb{R}^{L_{\rm{in}}\times N}}$, ${\bf{c}} \in {\mathbb{R}^{L_{\rm{in}}}}$, and $L_{\rm{in}}$ represents the number of linear inequalities. And we use $\preceq$ to denote component-wise inequality. The polyhedron is characterized by a set of linear inequalities, \textit{i.e.}, ${\bf{D}}{\bf{p}} \preceq {\bf{c}}$. Considering the cutting plane generation subproblem with polyhedron uncertainty set, which is required to solve the following problem,
\begin{align}
\label{eq:polyhedron_cutting}
& {\bf{p}}^{t+1}  = \mathop {\arg \max }\limits_{  {p_1},\dots,{p_N} }  \sum\nolimits_{j = 1}^N {({p_j}  - \mathop {{\rm{ }}p}\limits^ -  ){{f}_j}({{\boldsymbol{w}_j}})}\\
&{\rm{s.t.}}\quad \quad \quad  {\bf{D}}{\bf{p}} \preceq {\bf{c}} ,\sum\nolimits_{j = 1}^N {{p_j} }  = 1,  \nonumber\\
&{\rm{var.}}\quad \quad \quad \; {p_1},\dots,{p_N}.  \nonumber
\end{align}

The problem in Eq. (\ref{eq:polyhedron_cutting}) is an LP, which can be solved in polynomial time through interior point method \citep{boyd2004convex}. Specifically, the arithmetic complexity for interior point method to find the $\varepsilon'$-solution for LP is upper bounded by $\mathcal{O}((m\!+\!n)^{3/2}{n^{2}}\log (\frac{1}{{\varepsilon '}}))$ \citep{ben2011lectures}, where $m$ and $n$ are the number of inequality constraints and variables, respectively. In Eq. (\ref{eq:polyhedron_cutting}), $m\! =\! L_{\rm{in}}$ and $n\! = \!N$. The polyhedron uncertainty set cannot flexibly tradeoff between robustness with performance since it needs to adjust ${L_{\rm{in}} \! \times \! ({N\!+\!1})}$ parameters, \textit{i.e.}, ${\bf{D}}$ and ${\bf{c}}$.

\subsection{KL-Divergence Uncertainty Set}
The form of KL-divergence uncertainty set is given by,
\begin{align}
\label{eq:KL}
\boldsymbol{\mathcal{P}} =   \{ {\bf{p}}: \sum\limits_{j = 1}^N {{p_j}\log \frac{{{p_j}}}{{{q_j}}}}   \le {\beta}, {\bf{1}^ \top }{\bf{p}} = 1\} .
\end{align} 

Considering the cutting plane generation subproblem with KL-divergence uncertainty set, which is required to solve the following problem,
\begin{align}
\label{eq:KL_cutting}
& {\bf{p}}^{t+1}  = \mathop {\arg \max }\limits_{{p_1},\dots,{p_N} }  \sum\nolimits_{j = 1}^N {({p_j}  - \mathop {{\rm{ }}p}\limits^ -  ){{f}_j}({{\boldsymbol{w}_j}})}\\
&{\rm{s.t.}}\quad \;  \sum\limits_{j = 1}^N {{p_j}\log \frac{{{p_j}}}{{{q_j}}}}   \le {\beta},\sum\nolimits_{j = 1}^N {{p_j} }  = 1,  \nonumber\\
&{\rm{var.}}\quad \quad \quad \quad \; {p_1},\dots,{p_N}. \nonumber 
\end{align}

The above problem is a relative entropy programming (REP) \citep{chandrasekaran2017relative}. The arithmetic complexity for interior point method to find the $\varepsilon'$-solution for REP is bounded from above by $\mathcal{O}({(n)^{7/2}}|\log (\varepsilon ')|)$ \citep{potra1993quadratically} and $n=N$ in Eq. (\ref{eq:KL_cutting}). KL-divergence uncertainty set can flexibly tradeoff between robustness with performance by adjusting the parameter $\beta$. \textcolor{black}{ Additionally, KL-divergence uncertainty set supports statistical guarantees, ensuring robustness against moderate distributional shifts and often leading to improved out-of-sample performance \cite{namkoong2016stochastic}}.

\subsection{Earth-Mover (Wasserstein-1) Distance  Uncertainty Set}

KL-Divergence cannot deal with the prior distribution with zero elements \citep{qian2019robust} and have many drawbacks, \textit{e.g.}, asymmetry \citep{arjovsky2017wasserstein}. Earth-Mover (Wasserstein-1) distance becomes popular recently which can overcome the aforementioned drawbacks. The Earth-Mover (Wasserstein-1) distance uncertainty set can be expressed as,
\begin{align}
\label{eq:Wasserstein}
\boldsymbol{\mathcal{P}} \! = \!  \{ {\bf{p}}\!:\! \mathop {\min }\limits_{\gamma  \in \Pi ({\bf{p}},{\bf{q}})} \sum\limits_{i = 1}^N {\sum\limits_{j = 1}^N \! {\gamma ({x_i},{y_j})||{x_i} - {y_j}||} } \! \le \! {\beta},{\bf{1}^ \top }{\bf{p}} \!=\! 1\}  ,
\end{align}
where $\Pi ({\bf{p}},{\bf{q}})$ denotes the set of all joint distributions $\gamma $ whose marginal distributions are ${\bf{p}}$ and ${\bf{q}}$, respectively. And $x_i \! = \! i, \; i=1,\dots,N $, $y_j \!=\! j, \; j=1,\dots,N$. Intuitively, $\gamma ({x_i},{y_j})$ denotes the amount of ``mass'' that be moved from worker $i$ to  worker $j$ to transform the distribution ${\bf{p}}$ into the distribution ${\bf{q}}$. Thus, the Earth-Mover (Wasserstein-1) distance can be regarded as the ``cost'' of the optimal transport plan \citep{arjovsky2017wasserstein}. Considering the cutting plane generation subproblem with Earth-Mover (Wasserstein-1) distance uncertainty set, which is required to solve the following problem,
\begin{align}
\label{eq:Wasserstein_cutting}
& {\bf{p}}^{t+1}  = \mathop {\arg \max }\limits_{{p_1},\dots,{p_N}}  \sum\limits_{j = 1}^N {({p_j}  - \mathop {{\rm{ }}p}\limits^ -  ){{f}_j}({{\boldsymbol{w}_j}})}\\
&{\rm{s.t.}}\;   \sum\limits_{i = 1}^N {\sum\limits_{j = 1}^N {\gamma ({x_i},{y_j})||{x_i} - {y_j}||} } \! \le \! {\beta},\sum\limits_{j = 1}^N \! {{p_j} }  = 1,  \nonumber\\
&\quad \sum\limits_{j = 1}^N \! {\gamma ({x_i},{y_j})}  = {p_i},\forall i,  \sum\limits_{i = 1}^N \! {\gamma ({x_i},{y_j})}  = {q_j},\forall j, \nonumber\\
&{\rm{var.}}\quad  {p_1},\dots,{p_N}, {\gamma ({x_1},{y_1})},\dots,{\gamma ({x_N},{y_N})}. \nonumber 
\end{align}

The above problem is an LP. As mentioned above, the arithmetic complexity of finding the $\varepsilon'$-solution for problem (\ref{eq:Wasserstein_cutting}) through interior point method is bounded from above by $\mathcal{O}((m+n)^{3/2}{n^{2}}\log (\frac{1}{{\varepsilon '}}))$ \citep{ben2011lectures}, and  $m\! =\! 1$ and $n\! = \!N+N^2$ in this problem. It is seen that problem in Eq. (\ref{eq:Wasserstein_cutting}) is computationally more expensive to be solved since there is a quadratic
number of variables \citep{wong2019wasserstein}. The Earth-Mover (Wasserstein-1) distance uncertainty set can flexibly tradeoff between robustness with performance by adjusting one parameter, \textit{i.e.}, $\beta$. \textcolor{black}{Similar to KL-divergence uncertainty set, Wasserstein distance uncertainty set also has statistical benefits, such as robustness guarantees under distributional shifts \cite{gao2024wasserstein,liu2021stable} and dimension-free generalization bounds \cite{le2024universal}.}

{\color{black}
\subsection{How to Choose Different Uncertainty Set}}

Utilizing different uncertainty set results in different arithmetic complexity when solving the corresponding cutting plane generation subproblem. And the arithmetic complexity is also related to the number of workers $N$ (since $n$ in Table \ref{tab:uncertainty set comparison} is related to $N$). It is seen from \textcolor{black}{Figure 1 (a)} that the arithmetic complexity of KL-Divergence, ellipsoid, polyhedron and Wasserstein-1 distance uncertainty sets will increase significantly with the number of workers $N$. And the complexity of solving the cutting plane generation subproblem when utilizing Wasserstein-1 distance uncertainty set will increase significantly with the number of workers $N$, since there is a quadratic number of variables (\textit{i.e.}, $n=N+N^2$) in the cutting plane generation subproblem in Eq. (\ref{eq:Wasserstein_cutting}). And it is seen from \textcolor{black}{Figure 1 (b)} that, the complexity of utilizing ellipsoid uncertainty set will significantly increases when the number of constraints $m$ in cutting plane generation subproblem increases. Moreover, the complexity of utilizing   KL-Divergence, ellipsoid, polyhedron and Wasserstein-1 distance uncertainty sets will also increase quickly when $\varepsilon'$ decreases, which can be seen in \textcolor{black}{Figures 2 (a) and (b)}. As a result, from Figures 1 and 2, we can conclude that utilizing box and $CD$-norm uncertainty sets is more computationally efficient, especially when the distributed system is large (corresponding to a large $N$).

{\color{black}
In summary, the selection of uncertainty sets should be guided by the specific requirements in various scenarios. For instance,
\begin{enumerate}
    \item \textbf{Scenario 1: Limited master-side computational resources}. It is common for the master to have limited computational resources in Edge and IoT applications, as it is often deployed in resource-constrained edge data centers \cite{khan2023towards}. In this scenario, The box and $CD$-norm uncertainty sets can be utilized because the cutting plane generation subproblems are easier to solve when using box and $CD$-norm uncertainty sets.

    \item \textbf{Scenario 2: Low-Latency Requirement}. In real-time federated learning applications \cite{zhang2022scalable}, low-latency constitutes a critical performance metric. In this scenario, box and $CD$-norm uncertainty sets are preferred since the corresponding cutting-plane generation subproblems are simple linear programming, which can be effectively addressed by using merge sort.

    \item \textbf{Scenario 3: Flexibility}. For scenarios demanding flexible robustness adjustment, e.g., \cite{liu2024robust}, the ellipsoid, KL-divergence, Wasserstein-1 distance, and $CD$-norm uncertainty sets are preferable, as their robustness level can be flexibly regulated via a single parameter.

    \item \textbf{Scenario 4: Less conservatism}. In the scenario that less conservatism is required \cite{roos2020reducing}, meaning robustness should be ensured without compromising excessive optimality, the ellipsoid uncertainty set represents an appropriate selection, consistent with the discussion in \cite{bertsimas2004price}.

    \item \textbf{Scenario 5: Limited prior knowledge}. In scenarios with limited prior knowledge, particularly when the prior knowledge about $\widetilde{p}_j, \forall j$ is unavailable, uncertainty sets based on KL-divergence and Wasserstein-1 distance are preferable. Compared with the $CD$-norm uncertainty set, KL-divergence and Wasserstein-1 distance uncertainty sets do not require additional prior knowledge.

\end{enumerate}}

To facilitate the comparison of various uncertainty sets, the results of various uncertainty sets are summarized in Table \ref{tab:uncertainty set comparison}.}

\section{Convergence Analysis}\label{convergence analysis}

\noindent \textbf{Definition 1} {\rm{(Stationarity gap)}} Following \citep{xu2020unified,lu2020hybrid,xu2021zeroth}, the \textit{stationarity} \textit{gap} of our problem at $t^{{th}}$ iteration is defined as:
\begin{equation}
\nabla   G^t  =  \left[ \begin{array}{l}
   \{   \frac{1}{{{\alpha _{\boldsymbol{w}}^t}}}  (   {{\boldsymbol{w}}_j^t}  -  {\mathcal{P}_{{\boldsymbol{\mathcal{W}}}}} ( {{\boldsymbol{w}}_j^t}   -  {\alpha _{\boldsymbol{w}}^t} {\nabla _{{{\boldsymbol{w}}_j}}}      L_p^t))\vspace{0.3ex}\\
   
\, \frac{1}{{{\eta _{\boldsymbol{z}}^t}}}   (  {\boldsymbol{z}}^t    -    {\mathcal{P}_{{\boldsymbol{\mathcal{Z}}}}}(  {\boldsymbol{z}}^t    -    {\eta _{\boldsymbol{z}}^t} {\nabla _{\boldsymbol{z}}}  L_p^t)) \vspace{0.3ex}\\

\, \frac{1}{{{\eta _h^t}}}   (   h^t  -  {\mathcal{P}_{{\boldsymbol{\mathcal{H}}}}}(  h^t  -  {\eta _h^t} {\nabla _h} L_p^t))   \vspace{0.3ex}\\

\{  \frac{1}{{{\rho _1}}}( {\lambda _l^t}   -   {\mathcal{P}_{{\boldsymbol{\Lambda}}} }( {\lambda _l^t}   +    {\rho _1}{\nabla _{{\lambda _l}}}L_p^t))  \vspace{0.3ex}\\

\{   \frac{1}{{{\rho _2}}}( {{\boldsymbol{\phi}}_j^t}   -   {\mathcal{P}_{{\boldsymbol{\Phi}}}}( {{\boldsymbol{\phi}}_j^t}   +   {\rho _2}{\nabla _{{{\boldsymbol{\phi}}_j}}} L_p^t))   \end{array} \right],
\end{equation}
\noindent where $\nabla G^t$ and $L_p^t$ are respectively the simplified form of $\nabla G{\rm{(\{ }}{{\boldsymbol{w}}_j^t}{\rm{\} ,}}{\boldsymbol{z}^t},h^t,\{ {\lambda _l^t}\} ,\{ {{\boldsymbol{\phi}}_j^t}\} {\rm{)}}$ and ${L_p}{\rm{(\{ }}{{\boldsymbol{w}}_j^t}{\rm{\} ,}}{\boldsymbol{z}^t},h^t,\{ {\lambda _l^t}\} ,\{ {{\boldsymbol{\phi}}_j^t}\} {\rm{)}}$.

\noindent \textbf{Definition 2}
{\rm{($\varepsilon$-stationary point)}}  ${\rm{(\{ }}{{\boldsymbol{w}}_j^t}{\rm{\} }},{\boldsymbol{z}^t},h^t,\{ {\lambda _l^t}\} ,\{ {{\boldsymbol{\phi}}_j^t}\} {\rm{)}}$ is an $\varepsilon$-stationary point ($\varepsilon  \ge 0$) of a differentiable function ${L_p}$,  if $\,||\nabla G^t||^2 \le \varepsilon $.  $T(\varepsilon )$ is the first iteration index such that $||\nabla G^t||^2  \le  \varepsilon$, \textit{i.e.}, $T(\varepsilon )  =  \min \{ t \ |\; ||\nabla G^t||^2  \le  \varepsilon \}  $.

\noindent \textbf{Assumption 1} {\rm{(Smoothness/Gradient Lipschitz)}} $L_p$ has Lipschitz continuous gradients. We assume that there exists $L>0$ satisfying
\begin{equation}
\begin{array}{l}
\!||{\nabla }{L_p}( \{ {{\boldsymbol{w}}_j}\},\!{\boldsymbol{z}},\!h, \! \{ {\lambda _l}\}, \!\{ {{\boldsymbol{\phi}}_j}\} ) \! - \! {\nabla  }{L_p}( \{ {\hat{\boldsymbol{w}}_j}\},\!\hat{\boldsymbol{z}},\!\hat{h}, \! \{ {\hat{\lambda }_l}\}, \!\{ {\hat{\boldsymbol{\phi}}_j}\}  )|| \vspace{0.5ex}\\
 \le L||[ {{\boldsymbol{w}}_{\rm{cat}}} \! - \!  {{\hat{\boldsymbol{w}}_{\rm{cat}}}}   ;{\boldsymbol{z}} \! - \!{\hat{\boldsymbol{z}}}  ;h \! - \! \hat{h} ; {\boldsymbol{\lambda} _{\rm{cat}}} \! - \! {{\hat{\boldsymbol{\lambda}} _{\rm{cat}}}}  ; {{\boldsymbol{\phi}}_{\rm{cat}}} \! - \! {{\hat{\boldsymbol{\phi}}_{\rm{cat}}}}  ]||,
\end{array}
\end{equation}
where $[;]$ represents the concatenation. ${{\boldsymbol{w}}_{\rm{cat}}} \!-\! {{\hat{\boldsymbol{w}}_{\rm{cat}}}}   \!=\![{\boldsymbol{w}}_1 \! - \! {\hat{\boldsymbol{w}}_1};\cdots;{\boldsymbol{w}}_N \! - \! {\hat{\boldsymbol{w}}_N}] \! \in \! {\mathbb{R}^{pN}}$, ${\boldsymbol{\lambda} _{\rm{cat}}} \! - \! {{\hat{\boldsymbol{\lambda}} _{\rm{cat}}}} \!=\! [{\lambda_1} \! - \! {\hat{\lambda}_1};\cdots;{\lambda}_{|{{\bf{A}}^t}|} \! - \! {\hat{\lambda}_{|{{\bf{A}}^t}|}}]\! \in \! {\mathbb{R}^{|{{\bf{A}}^t}|}}$, ${{\boldsymbol{\phi}}_{\rm{cat}}}\!-\!{{\hat{\boldsymbol{\phi}}_{\rm{cat}}}} \!=\![{\boldsymbol{\phi}}_1 \! - \! {\hat{\boldsymbol{\phi}}_1};\cdots;{\boldsymbol{\phi}}_N \! - \! {\hat{\boldsymbol{\phi}}_N}]\! \in \! {\mathbb{R}^{pN}}$.

\noindent \textbf{Assumption 2} {\rm{(Boundedness)}} Before getting the $\varepsilon$-stationary point (\textit{i.e.}, $t\le T(\varepsilon )-1$), we assume variables in master satisfy that $||\boldsymbol{z}^{t+1} - \boldsymbol{z}^t|{|^2}+||h^{t + 1} - h^t|{|^2}+\sum\nolimits_l||{\lambda _l^{t + 1}} - {\lambda _l^t}|{|^2} \ge \vartheta $, where $\vartheta >0$ is a relative small constant. The change of the variables in master is upper bounded within $\tau$ iterations:
\begin{equation}
 \begin{array}{*{20}{l}}
{||\boldsymbol{z}^t - \boldsymbol{z}^{t - k}|{|^2} \! \le \! \tau{k_1}\vartheta}, \;\;
{||h^t - h^{t - k}|{|^2} \! \le \! \tau{k_1}\vartheta},\\
{\sum\nolimits_l||{\lambda _l^t} - {\lambda _l^{t - k}}|{|^2} \! \le \! \tau{k_1}\vartheta}, {\forall 1 \! \le \! k \! \le \! \tau },
\end{array}
\end{equation}
where $k_1 >0$ is a constant.


\noindent \textbf{Setting 1} {\rm{(Bounded $|{{\bf{A}}^t}|$)}} $|{{\bf{A}}^t}| \le M,{\rm{    }}\forall t$, \textit{i.e.}, an upper bound is set for the number of cutting planes.

\noindent \textbf{Setting 2} {\rm{(Setting of ${c_1^t}$, ${c_2^t}$)}}
${c_1^t}\! = \!\frac{1}{{{\rho _1}{(t+1)^{\frac{1}{6}}}}} \!\ge\! \underline{c}_1$ and ${c_2^t} \!=\! \frac{1}{{{\rho _2}{(t+1)^{\frac{1}{6}}}}} \!\ge\! \underline{c}_2$ are nonnegative non-increasing sequences, where $\underline{c}_1$ and $\underline{c}_2$ are positive constants.

\noindent \textbf{Theorem 1} {\rm{(Iteration complexity)}} Suppose Assumptions 1 and 2 hold. We set  
${\eta _{\boldsymbol{w}}^t} = {\eta _{\boldsymbol{z}}^t} = {\eta _h^t} = \frac{2}{{L + {\rho _1}|{{\bf{A}}^t}|{L^2} + {\rho _2}N{L^2} + 8(\frac{{|{{\bf{A}}^t}|\gamma {L^2}}}{{{\rho _1}({c_1^t})^2}} + \frac{{N\gamma {L^2}}}{{{\rho _2}({c_2^t})^2}})}}$ and $\underline{\eta _{\boldsymbol{w}}}=\frac{2}{{L + {\rho _1}M{L^2} + {\rho _2}N{L^2} + 8(\frac{{M\gamma {L^2}}}{{\rho _1}{\underline{c}_1}^2} + \frac{{N\gamma {L^2}}}{{\rho _2}{\underline{c}_2}^2})}}$. And we set constants  ${\rho _1} \!< \!\min \{\frac{ 2}{{L + 2c_1^0}}  ,\frac{1}{{15\tau{k_1}N{L^2}}}\} $ and $ \rho _2 \!\le\! \frac{2}{{L + 2c_2^0}} $, respectively. For a given $\varepsilon $, we have:
\begin{equation}
\begin{array}{l}
     T( \varepsilon ) \! \sim  \! \mathcal{O}(\max  \{ {(\frac{{4M\!{\sigma _1}^2}}{{{\rho _1}^2}}\! +\! \frac{{4N\!{\sigma _2}^2}}{{{\rho _2}^2}}\!)^3}\!\frac{1}{{{\varepsilon ^3}}}, \\

\qquad \qquad      
{(\frac{{4{{{(d_6+ \frac{{{\rho _2}(N  -  S){{L}^2}}}{2}\!)}\!}^2}\! (\mathop d\limits^ -  +  k_d(\tau \! - \! 1))  {d_5}}}{{{\varepsilon}}}\! +\! (T_1\!+\!\tau)^{\frac{1}{3}})^3}\}), 
\end{array}
\end{equation}
where ${\sigma _1}$, ${\sigma _2}$, $\gamma$, $k_d$, $\mathop d\limits^ -  $, ${d_5}$, ${d_6}$ and ${T_1}$ are constants. \textcolor{black}{Specifically, according to the setting of the proposed method, the dual variables $\lambda_l, {\boldsymbol{\phi}}_j$ are bounded, constants $\sigma_1$ and $\sigma_2$ represent the bounds for the dual variables, i.e., $||\lambda_l||\le \sigma_1, ||{\boldsymbol{\phi}}_j|| \le \sigma_2$. $\gamma$ is a constant that satisfies that $\gamma \ge \max\{2,  2 + \frac{{{\rho _1}{\rho _2}({c_1^0}{c_2^0}){^2}} (\max \{ \frac{1}{2},\frac{{3\tau{ k_1 }N{L^2}}}{2}\} - \frac{{{\rho _2}(N - S){L^2}}}{2})}{4L^2{{\rho _2}({c_2^0}){^2}} +4NL^2 {{\rho _1}({c_1^0}){^2}} }     \}$, and constant $d_6 = 4(\gamma  -  2){{L}^2}( M{\rho _1}  +  N{\rho _2})$. The constants $\mathop d\limits^ -$ and $k_d$, which depend on $M$, $N$, $\sigma_1$, $\sigma_2$, $\tau$, $\rho_1$, and $\rho_2$, are defined in detail in Eq. (95) of the Supplementary Materials. In addition, the constant $d_5$, which is related to $\tau$, $\rho_1$, $\rho_2$, $k_1$, $N$, $M$, $L$, and $\gamma$, is defined in Eq. (96) of the Supplementary Materials.}

\noindent \textbf{\textit{Proof:}} See Appendix A in Supplementary Materials.

There exists a wide array of works
regarding the convergence analysis of various algorithms for nonconvex/convex optimization problems involved in machine learning \citep{jin2020local,xu2021zeroth}.  Our analysis, however, differs from existing works in two aspects.  First, we solve the non-convex PD-DRO in an \emph{asynchronous} \emph{distributed} \emph{manner}. To our best knowledge, there are few works focusing on solving the DRO in a distributed manner. Compared to solving the non-convex PD-DRO in a centralized manner, solving it in an \emph{asynchronous} \emph{distributed} \emph{manner} poses significant challenges in algorithm design and convergence analysis. Secondly, we do not assume the inner problem can be solved nearly optimally for each outer iteration, which is numerically difficult to achieve in practice \citep{bertsekas1997nonlinear}. Instead, ASPIRE-EASE is \emph{single-loop} and involves simple gradient projection operations at each step.

{\color{black}
\noindent \textbf{Theorem 2} {\rm{(Communication complexity)}} Following \cite{backes2013asynchronous}, the communication complexity is defined as the cumulative count of bits exchanged between the involved parties, accounting for every transmitted message throughout the communication process. The overall communication complexity of the proposed method can be divided into the complexity at every iteration and the complexity of updating cutting planes. We have that the overall communication complexity of the proposed method is  $\mathcal{O}(\sum\nolimits_{t=1}^{T(\varepsilon)}{32S(2p+1+|{\bf{A}}^t|)} +\sum\nolimits_{t \in \mathcal{Q}} 32 N^2|{\bf{A}}^t|)$, where $\mathcal{Q}=\{k, 2k, \cdots, \lfloor \frac{T_1}{k} \rfloor k \}$.

\noindent \textbf{\textit{Proof:}} At every iteration, the active workers will update the local variables $\boldsymbol{w}_j^{t+1}$ and transmit $\boldsymbol{w}_j^{t+1}$ to the master. And master will update the variables and then broadcast the updated variables ${\boldsymbol{z}}^{t+1}$, $h^{t+1}$, $\{{\lambda _l^{t+1}}\}$ to active workers. Therefore, the communication complexity can be expressed as $C_1=32S(2p+1+|{\bf{A}}^t|)$.

Every $k$ iteration, the cutting planes will be updated in the master, the complexity of transmitting updated cutting plane parameters, i.e., ${\bf{A}}^t$, can be expressed as 
$C_2=\sum\nolimits_{t \in \mathcal{Q}} 32 N^2|{\bf{A}}^t|$, where $\mathcal{Q}=\{k, 2k, \cdots, \lfloor \frac{T_1}{k} \rfloor k \}$.

By combining with the iteration complexity $T(\varepsilon)$ at Theorem 1, we can obtain that the overall communication complexity of the proposed method is $\mathcal{O}(\sum\nolimits_{t=1}^{T(\varepsilon)} 32S(2p+1+|{\bf{A}}^t|)+\sum\nolimits_{t \in \mathcal{Q}} 32 N^2|{\bf{A}}^t|)$, where $\mathcal{Q}=\{k, 2k, \cdots, \lfloor \frac{T_1}{k} \rfloor k \}$.

}

\noindent \textcolor{black}{\textbf{Discussion.} In this work, since the objectives are non-convex, following the standard criterion in non-convex optimization, we utilize the \textit{first-order stationary point} to measure convergence \cite{qian2019robust}. Inspired by the fact that projected gradient descent and alternating gradient descent-ascent methods for non-convex optimization can achieve convergence guarantees, as demonstrated in \cite{liu2018zeroth,lin2020gradient}, the non-asymptotic convergence analysis is conducted by utilizing Cauchy-Schwarz
inequality, trigonometric inequality, projection
optimality condition, and Assumptions 1 and 2. The holistic analysis of the proofs of Theorems 1 and 2 is provided in Sec. VIII. In addition, compared with the related non-convex federated min-max learning works, the proposed ASPIRE-EASE is designed to handle the non-convex federated min-max optimization problem in an \textbf{asynchronous} manner. Conducting non-asymptotic convergence analysis for the proposed method is more challenging due to the varied update orders of individual workers, which introduce intricate interaction dynamics and significantly increase analytical complexity in the asynchronous algorithm. A comparison of the convergence rates between the proposed method and related non-convex federated min-max learning works regarding the convergence rate is shown in Table \ref{tab:example}.}

\renewcommand\arraystretch{1.2}
\renewcommand\tabcolsep{5pt}
\begin{table}[t]
\centering
\caption{\textcolor{black}{Convergence rate of federated learning works related to this work.}}
\begin{tabular}{c|c|c}
\toprule
Method & Synchronous min-max & Asynchronous min-max \\
\hline
\cite{sharma2022federated} & $\mathcal{O}(1/\varepsilon^{3.5})$ & NA \\
\hline
\cite{wu2023solving} & $\mathcal{O}(1/\varepsilon^{3})$ & NA \\
\hline
\cite{sharmafederated} & $\mathcal{O}(1/\varepsilon^{2})$ & NA \\
\hline
\cite{shen2024stochastic} & $\mathcal{O}(1/\varepsilon^{2})$ & NA \\
\hline
ASPIRE-EASE & NA & $\mathcal{O}(1/\varepsilon^3)$ \\
\bottomrule
\end{tabular}
\label{tab:example}
\end{table}

\begin{figure}[t]
\centering
\includegraphics[scale=0.5]{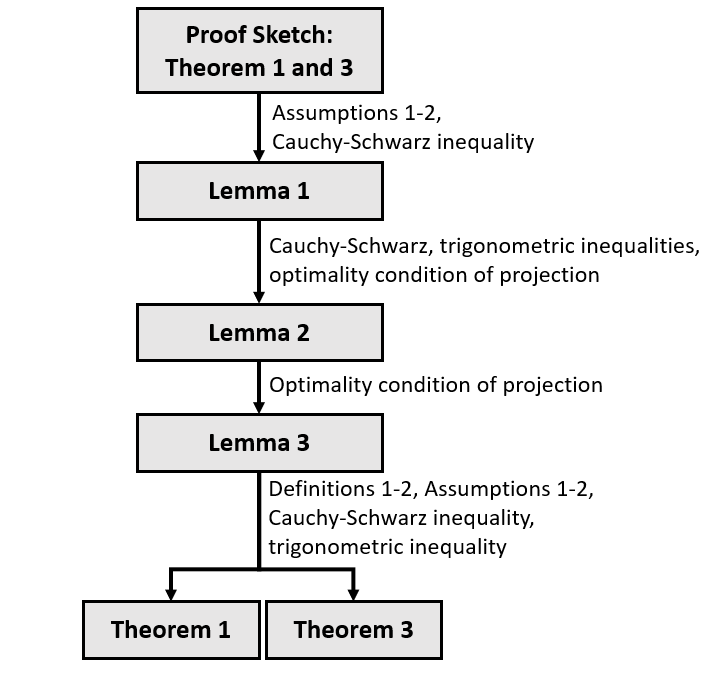}  
\caption{{\color{black}Key steps involved in the proofs of Theorems 1 and 3.}} 
\label{fig:holi1}
\end{figure}

\begin{table*}[t]
\caption{Unified complexity analysis for different uncertainty sets.}
\renewcommand\arraystretch{2.2}
\renewcommand\tabcolsep{15pt}
\label{tab:complexity}
\centering
{\scalebox{1}{\begin{tabular}{l|c} 
\toprule
Uncertainty Set &   Overall Arithmetic Complexity
\\
\hline
Box  & $\mathcal{O}\left(\frac{{\mathcal{K}_1}}{\varepsilon^3}   +  \left\lfloor {\frac{T_1}{k}} \right\rfloor n\log (n)\right)$ \\
\hline
Ellipsoid  & $\mathcal{O}\left(\frac{{\mathcal{K}_1}}{\varepsilon^3}   +  \left\lfloor {\frac{T_1}{k}} \right\rfloor (m+1)^{1/2}n({n^{2}+m+\sum\limits_{i = 1}^{m} {k_i^2} }) \log (\frac{1}{{\varepsilon '}})\right)$  \\
\hline
Polyhedron  & $\mathcal{O}\left(\frac{{\mathcal{K}_1}}{\varepsilon^3}   +  \left\lfloor {\frac{T_1}{k}} \right\rfloor (m+n)^{3/2}{n^{2}}\log (\frac{1}{{\varepsilon '}})\right)$  \\
\hline
KL-Divergence  & $\mathcal{O}\left(\frac{{\mathcal{K}_1}}{\varepsilon^3}   +  \left\lfloor {\frac{T_1}{k}} \right\rfloor{(n)^{7/2}}|\log (\varepsilon ')|\right)$ \\
\hline
Wasserstein-1 Distance   & $\mathcal{O}\left(\frac{{\mathcal{K}_1}}{\varepsilon^3}   +  \left\lfloor {\frac{T_1}{k}} \right\rfloor(m+n)^{3/2}{n^{2}}\log (\frac{1}{{\varepsilon '}})\right)$ \\
\hline
CD-norm    & $\mathcal{O}\left(\frac{{\mathcal{K}_1}}{\varepsilon^3}   +  \left\lfloor {\frac{T_1}{k}} \right\rfloor n\log (n)\right)$ \\
\bottomrule
\end{tabular}}}
\\\footnotesize{{\color{black}${\mathcal{K}_1}$ denotes the arithmetic complexity of gradient projections from Eq. (\ref{eq:15}) to Eq. (\ref{eq:y_update}).}}
\end{table*}

{\color{black}
\section{Unified Complexity Analysis}
\label{Unified Complexity Analysis}
In this section, we extend the proposed ASPIRE-EASE algorithm to different uncertainty sets and make a unified analysis regarding its arithmetic complexity. Recall that according to Theorem 1, the iteration complexity of the proposed method is upper bounded by $\mathcal{O}(\frac{1}{\varepsilon^3})$. In the first $T_1$ iterations, the cutting planes will be updated by solving the cutting plane generation subproblem for every $k$ iterations. Notice that the overall arithmetic complexity of the proposed algorithm is dominated by complexity of the gradient projection operations (from Eq. (\ref{eq:15}) to Eq. (\ref{eq:y_update})) and solving the cutting plane generation subproblem. Consequently, we can obtain the arithmetic complexity of the proposed algorithm with different uncertainty sets are summarized in Table \ref{tab:complexity}.}

{\color{black} 
\section{ASPIRE-ADP}
\label{ASPIRE-ADP}
We next propose ASPIRE-ADP, i.e., ASPIRE-EASE with an adaptive NAW. It is seen that such an adaptive technique can effectively accelerate the converge of the proposed algorithm. The setting of $S$, \textit{i.e.}, the number of active workers, controls the level of asynchrony of ASPIRE-EASE and has a direct impact on the training efficiency. For example, if we set $S=1$, the proposed algorithm is fully asynchronous. Likewise, if we let $S=N$, we obtain a fully synchronous distributed algorithm. Consequently, choosing a proper $S$ is crucial for the training efficiency. For instance, if all workers have almost the same capacities of computation and communication, \textit{i.e.}, each worker has roughly the same computation and communication 
delay, a fully synchronous algorithm, i.e., $S=N$ requires less training rounds) and easier to implement \citep{assran2020advances}. And in the other case, when there are stragglers in the distributed system, the algorithm will be more efficient if we set $S<N$ instead of $S=N$ \citep{zhang2014asynchronous}.

\begin{figure}[t]
\centering
\includegraphics[scale=0.45]{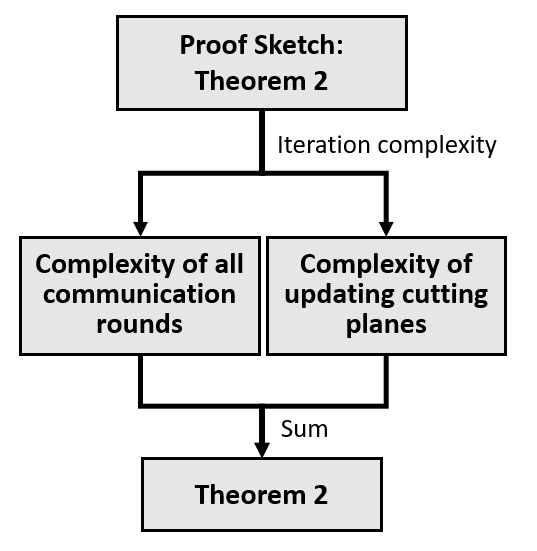}  
\caption{{\color{black}Key steps involved in the proofs of Theorem 2.}} 
\label{fig:holi2}
\vspace{-4mm}
\end{figure}

Thus, setting a proper $S$ in the distributed system is crucial. Nevertheless, the delay of some workers may change abruptly in the process of training. As a result, a fixed NAW may not be the optimal choice, adaptive NAW is more preferred. For example, \citep{su2022gba} has pointed out that switching the mode between synchronous and asynchronous training will enhance the training efficiency. We therefore propose ASPIRE-ADP. Specifically, the $S$ will be updated in master based on the estimated delay information of each worker as follows,
\begin{equation}
\label{eq:7_30_36}
S = \left\{ \begin{array}{l}
s,\;{\rm{ if  }}\quad \max \{ \mathcal{T}_j \}  \le {\beta _1}\\
N,\;{\rm{ if  }}\quad \max \{ \mathcal{T}_j\}   > {\beta _1}
\end{array} \right.,
\end{equation}
where $1\le s<N$ is an integer. \textcolor{black}{$\mathcal{T}_j$ denotes the estimated delay of worker $j, \forall j$, and master can obtain $\mathcal{T}_j$ based on the communication time interval with worker $j$. 
Specifically, let $\underline{e_j}$ denote the time of the last communication between worker $j$ and the master, $e$ denote the current time, and $\overline{e_j}$ denote the time of the next communication between worker $j$ and the master ($\underline{e_j} \le e \le \overline{e_j} $). Thus, the $\mathcal{T}_j$ can be obtained on the master as follows.
\begin{equation}
\mathcal{T}_j = \left\{ \begin{array}{l}
\max\{e-\underline{e_j}, \mathcal{T}_j\},\;{\rm{ when  }}\quad e<\overline{e_j}\\
\overline{e_j}-\underline{e_j},\;{\rm{ when  }}\quad e=\overline{e_j}
\end{array} \right. , \forall j.
\end{equation} }

\textcolor{black}{In Eq. (\ref{eq:7_30_36}), $\beta_1$ serves as a threshold to determine the presence of stragglers in a distributed system, thereby deciding whether a mode switch is necessary. $\beta _1$ can be flexibly selected based on the requirements of the distributed system. For example, 1) a fixed time threshold (e.g., 10 seconds) can be employed as $\beta_1$, under the rationale that any worker failing to communicate with the master node for an extended period can be considered a straggler \cite{xu2018chronos}. 2) $\beta_1$ can also be selected based on the average delay ${\rm{AVE}} \{ \mathcal{T}_j\}$ in the distributed system. As discussed in \cite{ouyang2016straggler},  a worker is identified as a straggler if its delay reaches twice the average delay, thus we can select $\beta _1 = 2 {\rm{AVE}} \{ \mathcal{T}_j\}$. }

We give an example to show how ASPIRE-ADP works, which can be seen in Figure \ref{fig:adaptive}. In the experiment, we assume that there are three workers in a distributed system. In time intervals 0$\sim$900s and 1300$\sim$2000s, all workers have the similar time delay. In the time interval 900$\sim$1300s, there is a straggler (\textit{i.e.}, worker 3), which leads to larger delay than the other two workers. As shown in Figure \ref{fig:adaptive}, ASPIRE-ADP can adjust NAW according to the estimated delay information.

Then, we analyze the iteration complexity of ASPIRE-ADP, as shown in Theorem 3. With an adaptive NAW, from $T_1+\tau$ iteration to $T( \varepsilon )$ iteration, we assume that the number of iterations when $S=s$ is $\beta_2(T( \varepsilon )-T_1-\tau+1)$ and the number of iterations when $S=N$ is $(1-\beta_2)(T( \varepsilon )-T_1-\tau+1)$, where $ 0\le \beta_2 \le 1$.

\begin{figure}[t]
\centering
\subfigure[Three workers with different delay.]
{\begin{minipage}{6.5cm}
    \includegraphics[scale=0.34]{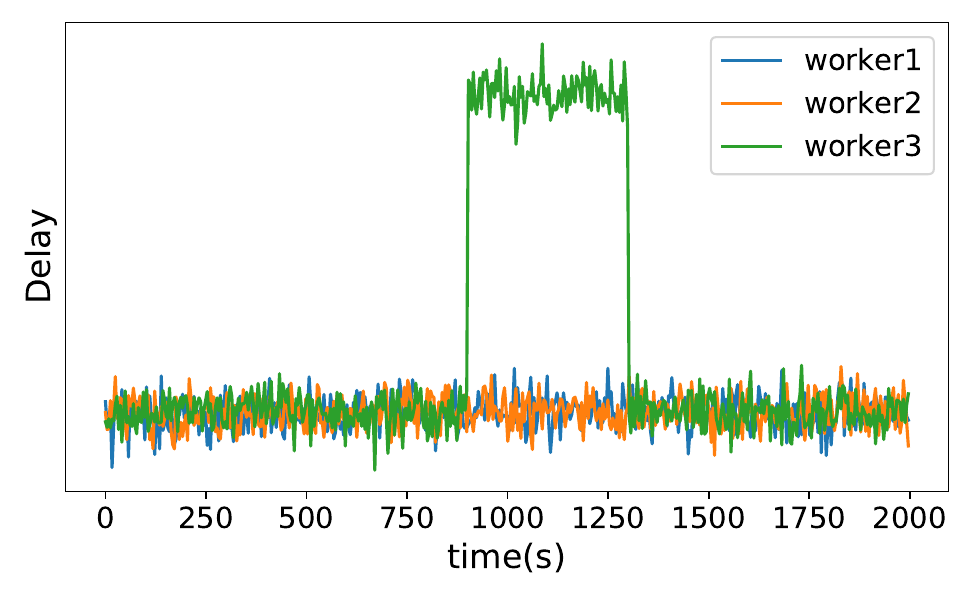}  
\end{minipage}}

\subfigure[\textcolor{black}{$S$ adaptively changes.}]
{\begin{minipage}{6.85cm}
    \includegraphics[scale=0.35]{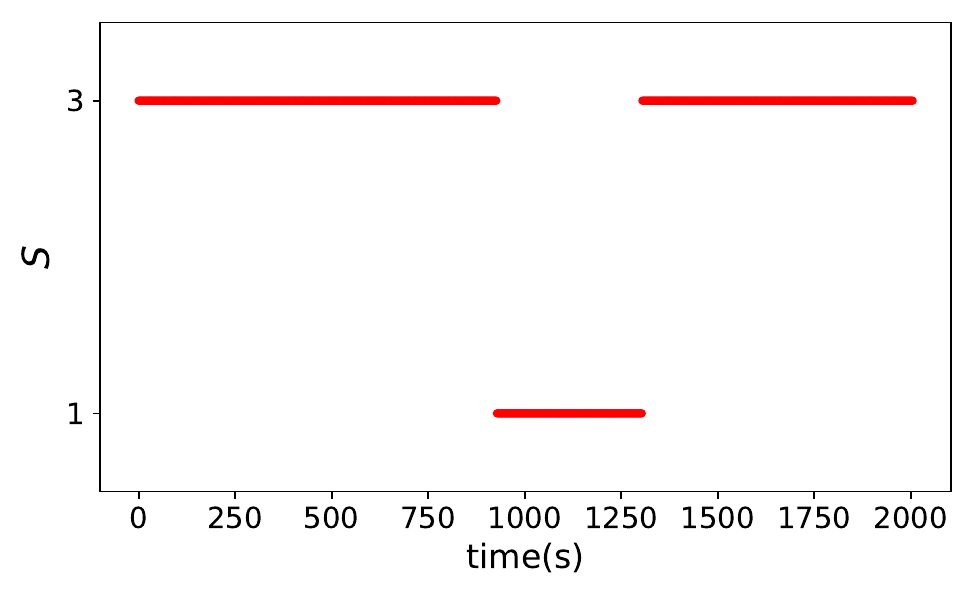}  
\end{minipage}}

\subfigure[\textcolor{black}{ASPIRE-ADP vs ASPIRE-EASE on Fashion MNIST dataset.}]
{\begin{minipage}{7.4cm}
    \includegraphics[scale=0.365]{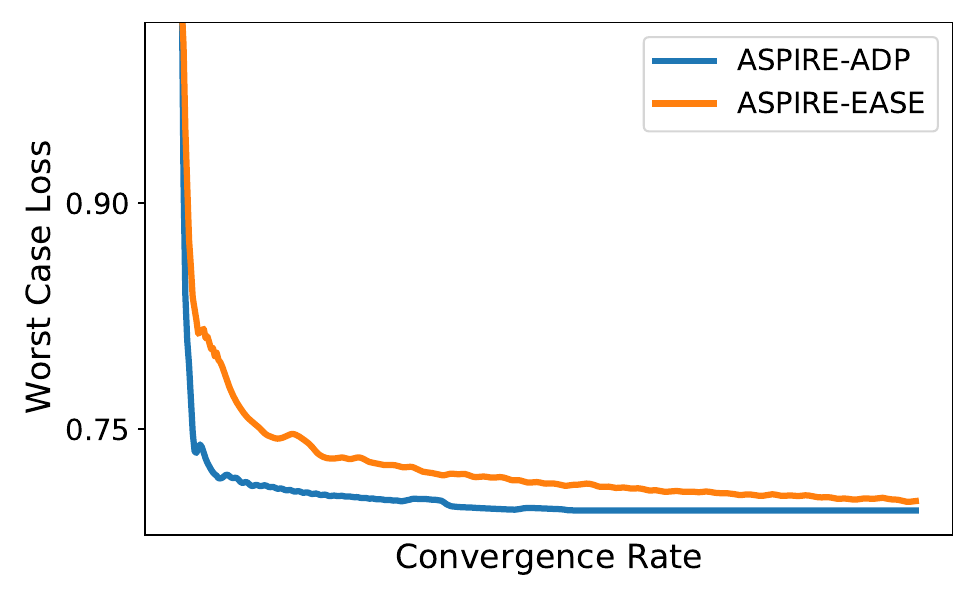}  
\end{minipage}}

\caption{{\color{black}The number of active workers $S$ can adjust adaptively based on the estimated delay information of each worker. In (a), there are three workers with different delays. Worker 3 is a straggler during time 900$\sim$1300s. In (b), $S$ changes accordingly based on the estimated delay information. In (c), we compare the iteration complexity of ASPIRE-ADP and ASPIRE-EASE.}} 
\label{fig:adaptive}
\end{figure}

{\color{black}\noindent \textbf{Theorem 3} {\rm{(Iteration complexity)}} Suppose Assumptions 1 and 2 hold. We set  
${\eta _{\boldsymbol{w}}^t} = {\eta _{\boldsymbol{z}}^t} = {\eta _h^t} = \frac{2}{{L + {\rho _1}|{{\bf{A}}^t}|{L^2} + {\rho _2}N{L^2} + 8(\frac{{|{{\bf{A}}^t}|\gamma {L^2}}}{{{\rho _1}({c_1^t})^2}} + \frac{{N\gamma {L^2}}}{{{\rho _2}({c_2^t})^2}})}}$ and $\underline{\eta _{\boldsymbol{w}}}=\frac{2}{{L + {\rho _1}M{L^2} + {\rho _2}N{L^2} + 8(\frac{{M\gamma {L^2}}}{{\rho _1}{\underline{c}_1}^2} + \frac{{N\gamma {L^2}}}{{\rho _2}{\underline{c}_2}^2})}}$. And we set constants  ${\rho _1} \!< \!\min \{\frac{ 2}{{L + 2c_1^0}}  ,\frac{1}{{15\tau{k_1}N{L^2}}}\} $ and $ \rho _2 \!\le\! \frac{2}{{L + 2c_2^0}} $, respectively. For a given $\varepsilon $, we have:
\begin{equation}
\begin{array}{l}
T( \varepsilon ) \! \sim  \! \mathcal{O}(\max  \{ {(\frac{{4M\!{\sigma _1}^2}}{{{\rho _1}^2}}\! +\! \frac{{4N\!{\sigma _2}^2}}{{{\rho _2}^2}}\!)^3}\!\frac{1}{{{\varepsilon ^3}}},
\\
\qquad \qquad
{(\frac{{4(\mathop d\limits^ -  + k_d(\tau -1))   {d_5}}}{{(\frac{\beta_2}{{(d_6 + \frac{{{\rho _2}(N - s){L^2}}}{2})^2}} + \frac{1-\beta_2}{{d_6}^2}){\varepsilon}}} + (T_1+\tau)^{\frac{1}{3}})^3}\}),
\end{array}
\end{equation}
where ${\sigma _1}$, ${\sigma _2}$, $\gamma$, $k_d$, $\mathop d\limits^ -  $, ${d_5}$, ${d_6}$ and ${T_1}$ are constants.  It is seen from Theorem 3 that ASPIRE-ADP can effectively improve the training efficiency and reduce the iteration complexity compared with ASPIRE-EASE.

\noindent \textbf{\textit{Proof:}} See Appendix B in Supplementary Materials.
}

}

\noindent \textcolor{black}{\textit{\textbf{Holistic analysis of the proofs of Theorems 1, 2, and 3.}} The proofs of iteration complexity (i.e., Theorem 1 and 3) are built upon Lemmas 1, 2, and 3, and can be summarized in four key steps. 1) First, by leveraging Assumptions 1 and 2 together with the Cauchy-Schwarz inequality, we can get Lemma 1. 2) Next, based on Lemma 1, Lemma 2 is derived by further incorporating Assumptions 1 and 2, the optimality condition of projection, the Cauchy–Schwarz inequality, and trigonometric inequality. 3) Then, by utilizing Lemma 2 and the projection optimality condition, Lemma 3 can be obtained. 4) Finally, by combining Lemma 3 with the definitions of stationarity gap and $\varepsilon$-stationary point, Assumptions 1 and 2, trigonometric inequality, and Cauchy-Schwarz  inequality, we can get the iteration complexity. The proofs of communication complexity (i.e., Theorem 2) are based on two steps. 1) The analyses of the communication complexity of all communication rounds based on iteration complexity, and 2) the analyses of the communication complexity of updating cutting planes. To enhance readability and facilitate a better understanding of the derivation process, we illustrate the key steps involved in the proofs of Theorems 1, 2, and 3 in Figures \ref{fig:holi1} and \ref{fig:holi2}.}

\noindent \textcolor{black}{\textbf{Discussion.} It is seen from Theorems 1 and 3 that ASPIRE-ADP demonstrates lower iteration complexity than ASPIRE-EASE, leading to improved computational efficiency. Detailed discussions about the specific reduced iteration complexity and costs across all communication rounds are provided in Appendix D.}

\section{Experiment} \label{experiment}

In this section, we conduct experiments on four real-world datasets to assess the performance of the proposed method. Specifically, we evaluate the robustness against data heterogeneity, robustness against malicious attacks, and efficiency of the proposed method. Ablation study is also carried out to demonstrate the excellent performance of ASPIRE-EASE.

\subsection{Datasets and Baseline Methods}\label{dataset}

We compare the proposed ASPIRE-EASE with baseline methods based on SHL~\citep{gjoreski2018university}, Person Activity~\citep{kaluvza2010agent}, Single Chest-Mounted Accelerometer (SM-AC)~\citep{casale2012personalization} and Fashion MNIST~\citep{xiao2017fashion} datasets. The baseline methods include ${\rm{Ind}}_j$ (learning the model from an individual worker $j$), ${\rm{Mix}}{\rm{_{Even}}}$ (learning the model from all workers with even weights using ASPIRE),  FedAvg \citep{mcmahan2017communication}, AFL \citep{mohri2019agnostic} and DRFA-Prox \citep{deng2021distributionally}. The detailed descriptions of datasets and baselines are given in Appendix C in Supplementary Materials.

In our empirical studies, since the downstream tasks are multi-class classification, the cross entropy loss is used on each worker (\textit{i.e.}, ${\mathcal{L}_j}( \cdot ),\forall j$). For SHL, Person Activity, and SM-AC datasets, we adopt the deep multilayer perceptron~\citep{wang2017time} as the base model. And we use the same logistic regression model as in \citep{mohri2019agnostic,deng2021distributionally} for Fashion MNIST dataset. The base models are trained with SGD. Following related works in this direction~\citep{qian2019robust,mohri2019agnostic,deng2021distributionally}, worst case performance are reported for the comparison of robustness. Specifically, we use {${\bf{Acc}}_{w}$} and {${\bf{Loss}}_{w}$} to represent the worst case test accuracy and training loss (\textit{i.e.}, the test accuracy and training loss on the worker with worst performance), respectively. We also report the standard deviation ${\bf{Std}}$ of $[\rm{Acc}_1,\cdots , \rm{Acc}_N]$ (the test accuracy on every worker). In the experiment, $S$ is set as 1, that means the master will make an update once it receives a message.  Each experiment is repeated 10 times, both mean and standard deviations are reported.

\renewcommand\arraystretch{1.3}
\renewcommand\tabcolsep{4pt}
\begin{table*}[t]
\caption{Performance comparisons based on {${\bf{Acc}}_{w}$} (\%) $\uparrow$, ${\bf{Loss}}_{w}$ $\downarrow$ and ${\bf{Std}}$ $\downarrow$  ($\uparrow$ and $\downarrow$ respectively denote higher scores represent better performance and lower scores represent better performance). The boldfaced digits represent the best results, ``$-$'' represents not available.}
\centering
\label{tab:III}
\scalebox{0.9}{
\begin{tabular}{l|cccccccccccc}
\toprule
\multirow{2}{*}{Model}&  
    \multicolumn{3}{c}{SHL}&\multicolumn{3}{c}{Person Activity}&\multicolumn{3}{c}{SC-MA}&\multicolumn{3}{c}{Fashion MNIST}\cr  
    \cmidrule(lr){2-4} \cmidrule(lr){5-7} 
    \cmidrule(lr){8-10}
    \cmidrule(lr){11-13}
    & ${\bf{Acc}}_{w}$$\uparrow$  & ${\bf{Loss}}_{w}$$\downarrow$  & {\bf{Std}}$\downarrow$  & ${\bf{Acc}}_{w}$$\uparrow$ & ${\bf{Loss}}_{w}$ $\downarrow$ & {\bf{Std}}$\downarrow$  & ${\bf{Acc}}_{w}$ $\uparrow$  & ${\bf{Loss}}_{w}$ $\downarrow$ & {\bf{Std}}$\downarrow$   & ${\bf{Acc}}_{w}$ $\uparrow$  & ${\bf{Loss}}_{w}$ $\downarrow$ & {\bf{Std}}$\downarrow$ \\ \hline
${\rm{max}}\{{\rm{Ind}}_j\}$  & 19.06±0.65 &  $-$  &  29.1   & 49.38±0.08  & $-$  &  8.32  & 22.56±0.78  & $-$  &   17.5  & $-$   &  $-$ &  $-$   \\
${\rm{Mix}}{\rm{_{Even}}}$   & 69.87±3.10  & 0.806±0.018  &  4.81    & 56.31±0.69  & 1.165±0.017  & 3.00     & 49.81±0.21  &  1.424±0.024 &   6.99   & 66.80±0.18   & 0.784±0.003  &   10.1   \\ 
FedAvg \citep{mcmahan2017communication}  & 69.96±3.07  & 0.802±0.023  &   5.21   & 56.28±0.63  & 1.154±0.019  &   3.13   & 49.53±0.96  & 1.441±0.015  &   7.17  & 66.58±0.39   & 0.781±0.002  &   10.2   \\ 
AFL  \citep{mohri2019agnostic}    & 78.11±1.99   & 0.582±0.021  &   1.87   & 58.39±0.37   &  1.081±0.014 &    0.99   & 54.56±0.79   & 1.172±0.018  &   3.50   & 77.32±0.15   & 0.703±0.001  &   1.86   \\
DRFA-Prox \citep{deng2021distributionally}    & 78.34±1.46   & 0.532±0.034  &   1.85    & 58.62±0.16  & 1.096±0.037  &   1.26   & 54.61±0.76  &  1.151±0.039 &  4.69   & 77.95±0.51 &  0.702±0.007 &   1.34   \\ \hline
ASPIRE-EASE     & \textbf{79.16±1.13}  & \textbf{0.515±0.019}  &  \textbf{1.02}    & 59.43±0.44  &   1.053±0.010 &  0.82    &  56.31±0.29 &  1.127±0.021 &  \textbf{3.16}  &     \textbf{78.82±0.07}  &  \textbf{0.696±0.004} &   \textbf{1.01}   \cr
ASPIRE-EASE$_{\rm{per}}$     &  78.94±1.27    &  0.521±0.023 &  1.36      & \textbf{59.54±0.21}  &   \textbf{1.051±0.016} &  \textbf{0.79}    &  \textbf{56.71±0.16} &  \textbf{1.119±0.028} &  3.48  &     78.73±0.06  &  0.698±0.006 &   1.09 \cr
\bottomrule  
\end{tabular}}
\end{table*}

\begin{figure}[t]
\subfigure[Person Activity]
{\begin{minipage}{4.2cm}
    \includegraphics[scale=0.285]{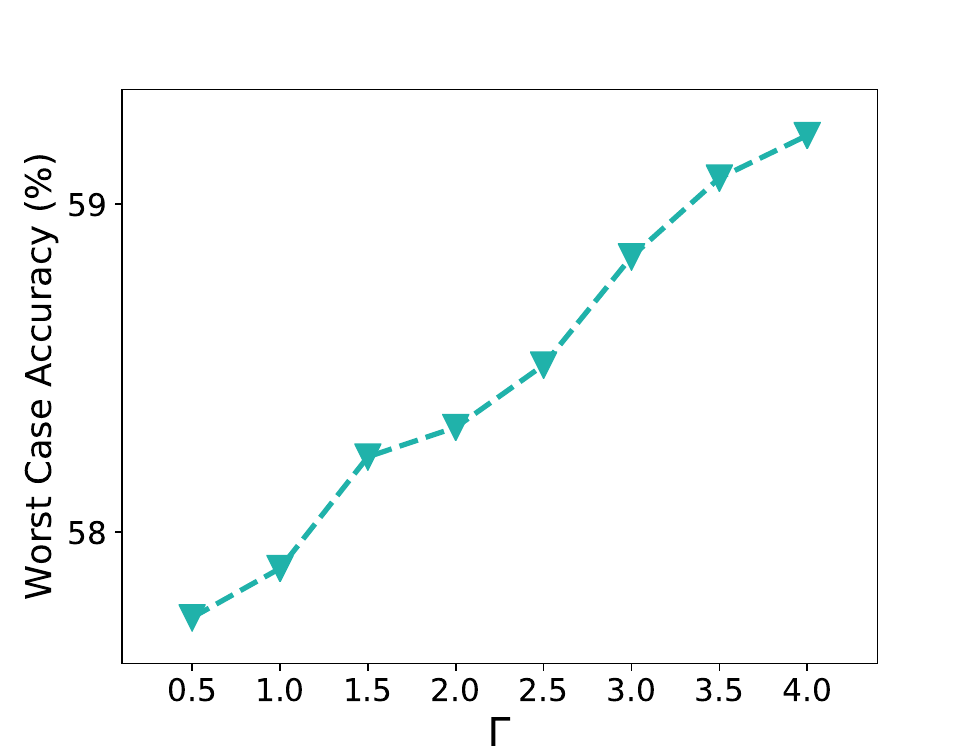}  \vspace{0.01mm}
\end{minipage}}
\subfigure[SC-MA] 
{\begin{minipage}{4.2cm}
    \includegraphics[scale=0.285]{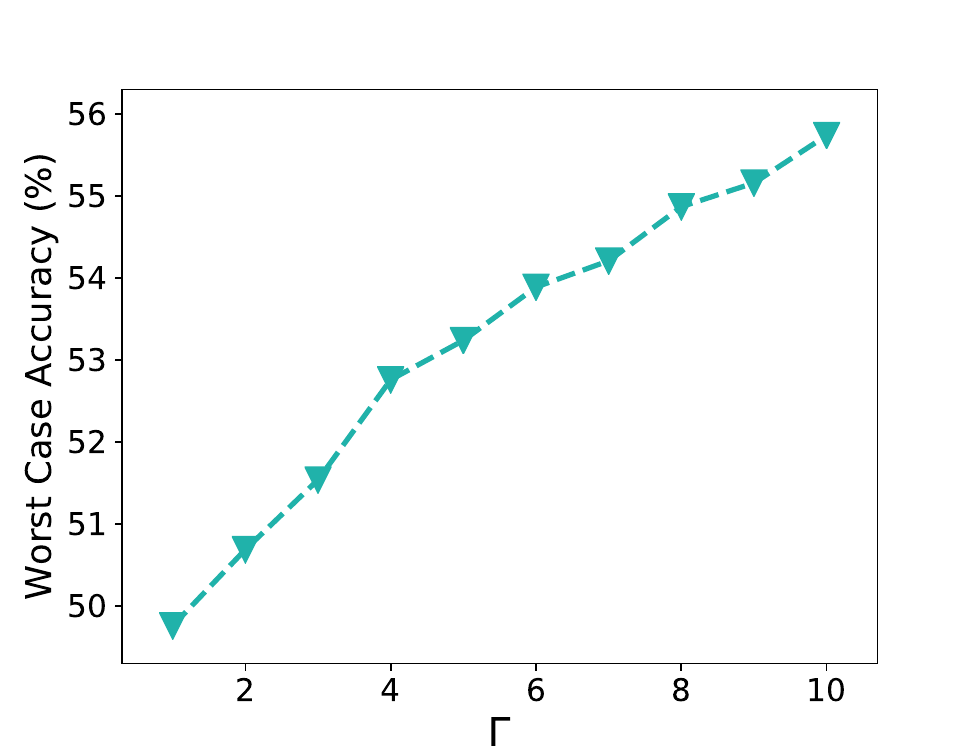}  \vspace{0.01mm}
\end{minipage}}

\caption{$\Gamma$ control the degree of robustness (worst case performance in the problem) on  (a) Person Activity, (b) SC-MA  datasets.} 
\label{fig:gamma1}
\end{figure}

\begin{figure}[t]
\subfigure[SHL]
{\begin{minipage}{4.2cm}
    \includegraphics[scale=0.285]
    {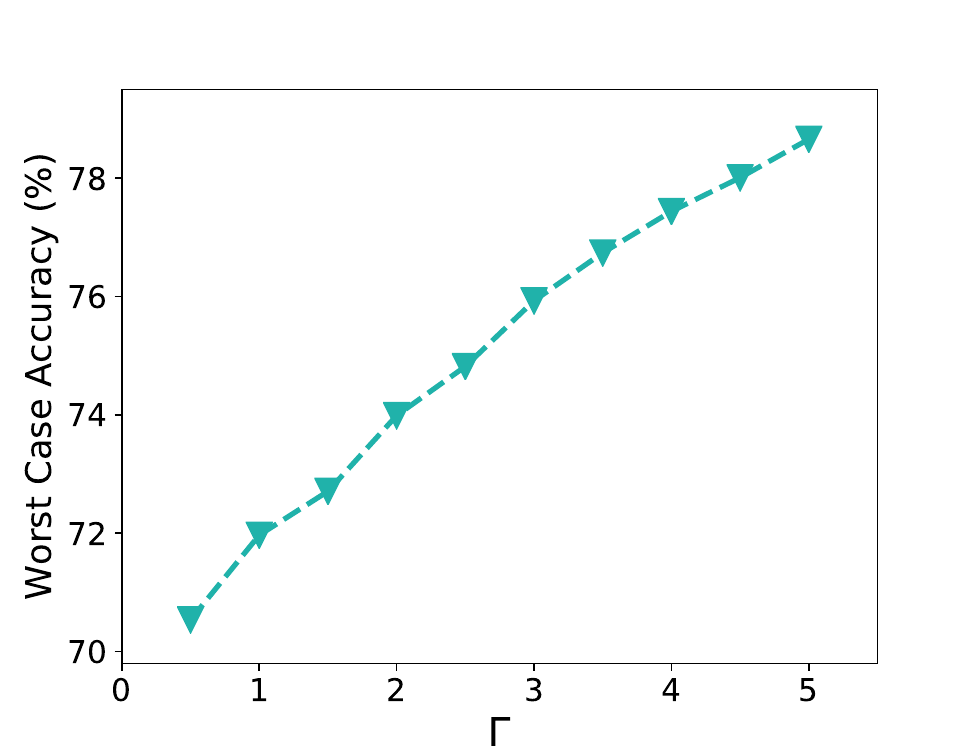} \vspace{0.01mm}
\end{minipage}}
\subfigure[Fashion MNIST] 
{\begin{minipage}{4.2cm}
    \includegraphics[scale=0.285]{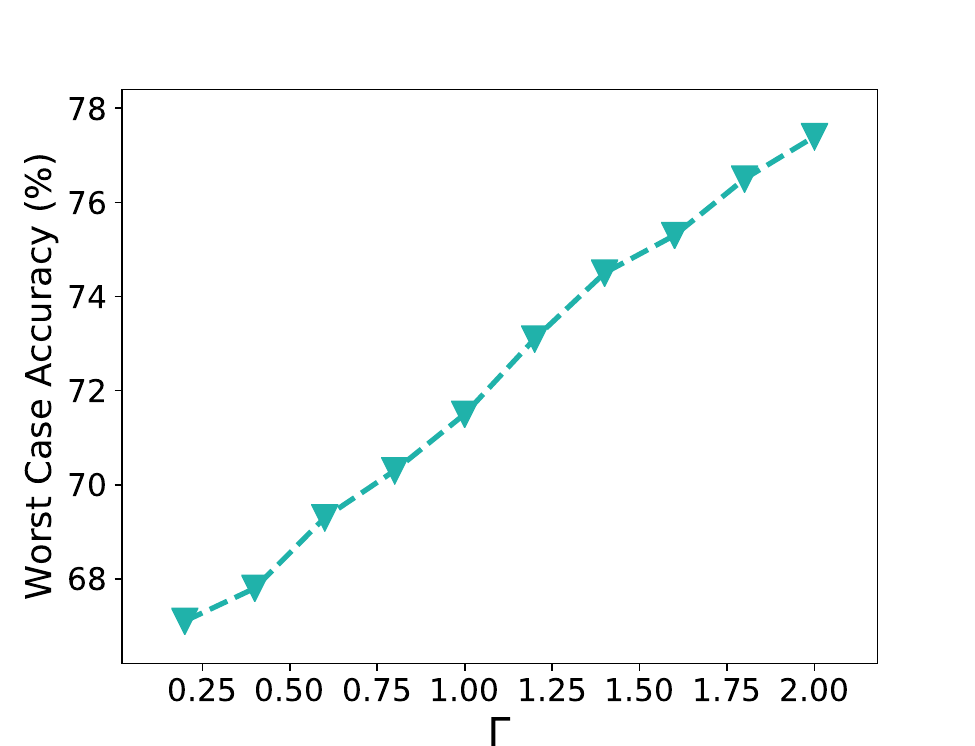} \vspace{0.01mm}
\end{minipage}}

\caption{$\Gamma$ control the degree of robustness (worst case performance in the problem) on  (a) SHL, (b) Fashion MNIST  datasets.} 
\label{fig:gamma2}
\end{figure}

\subsection{Results}

\subsubsection{Robustness against Data Heterogeneity} 
We first assess the robustness of the proposed ASPIRE-EASE by comparing it with baseline methods when data are heterogeneously distributed across different workers. Specifically, we compare the ${\bf{Acc}}_{w}$, ${\bf{Loss}}_{w}$ and ${\bf{Std}}$ of different methods on all datasets. The performance comparison results are shown in Table~\ref{tab:III}.  In this table, we can observe that ${\rm{max}}\{{\rm{Ind}}_j\}$, which represents the best performance of individual training over all workers, exhibits the worst robustness on SHL, Person Activity, and SC-MA. This is because individual training (${\rm{max}}\{{\rm{Ind}}_j\}$) only learns from the data in its local worker and cannot generalize to other workers due to different data distributions. Note that ${\rm{max}}\{{\rm{Ind}}_j\}$ is unavailable for Fashion MNIST since each worker only contains one class of data and cross entropy loss cannot be used in this case. ${\rm{max}}\{{\rm{Ind}}_j\}$ also does not have ${\bf{Loss}}_{w}$, since ${\rm{Ind}}_j$ is trained only on individual worker $j$. The FedAvg and ${\rm{Mix}}{\rm{_{Even}}}$ exhibit better performance than ${\rm{max}}\{{\rm{Ind}}_j\}$ since they consider the data from all workers. Nevertheless, FedAvg and  ${\rm{Mix}}{\rm{_{Even}}}$ only assign the fixed weight for each worker. 
AFL is more robust than FedAvg and  ${\rm{Mix}}{\rm{_{Even}}}$ since it not only utilizes the data from all workers but also considers optimizing the weight of each worker. DRFA-Prox outperforms AFL since it also considers the prior distribution and regards it as a regularizer in the objective function. Finally, we can observe that the proposed ASPIRE-EASE shows excellent robustness, which can be attributed to two factors: 1) ASPIRE-EASE considers data from all workers and can optimize the weight of each worker;
2) compared with DRFA-Prox which uses prior distribution as a regularizer, the prior distribution is incorporated within the constraint in our formulation (Eq. \ref{eq:3}), which can be leveraged more effectively. And it is seen that ASPIRE-EASE can perform periodic communication since ASPIRE-EASE$_{\rm{per}}$, which represents ASPIRE-EASE with periodic communication, also has excellent performance.

\begin{figure}[t]
\subfigure[Person Activity]
{\begin{minipage}{4.2cm}
    \includegraphics[scale=0.285]{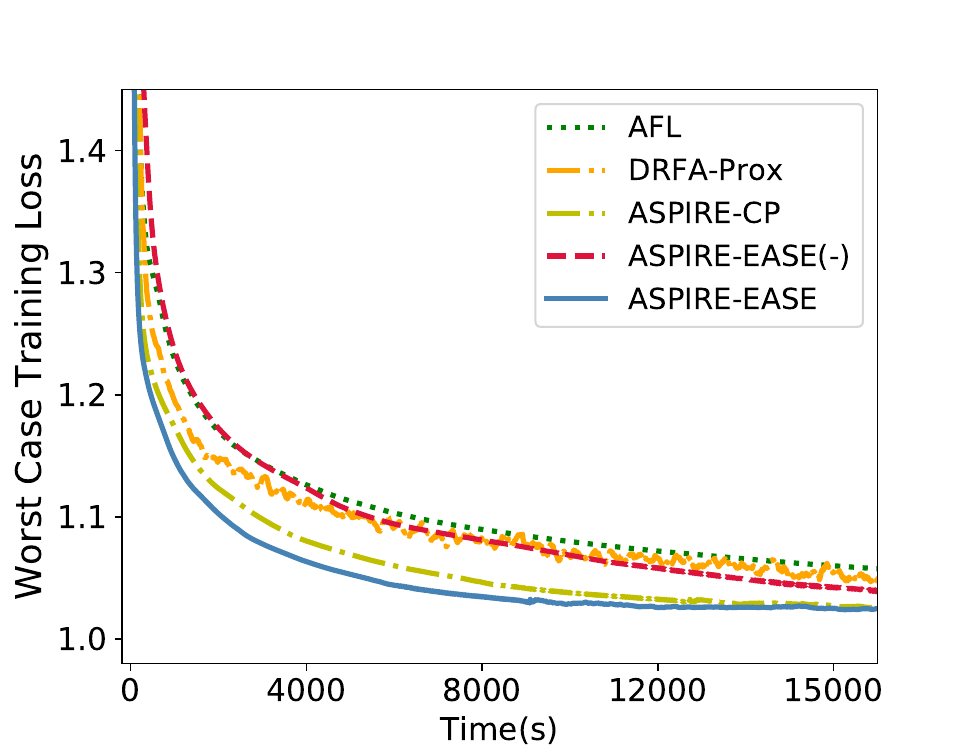}  \vspace{0.01mm}
\end{minipage}}
\subfigure[SC-MA] 
{\begin{minipage}{4.2cm}
    \includegraphics[scale=0.285]{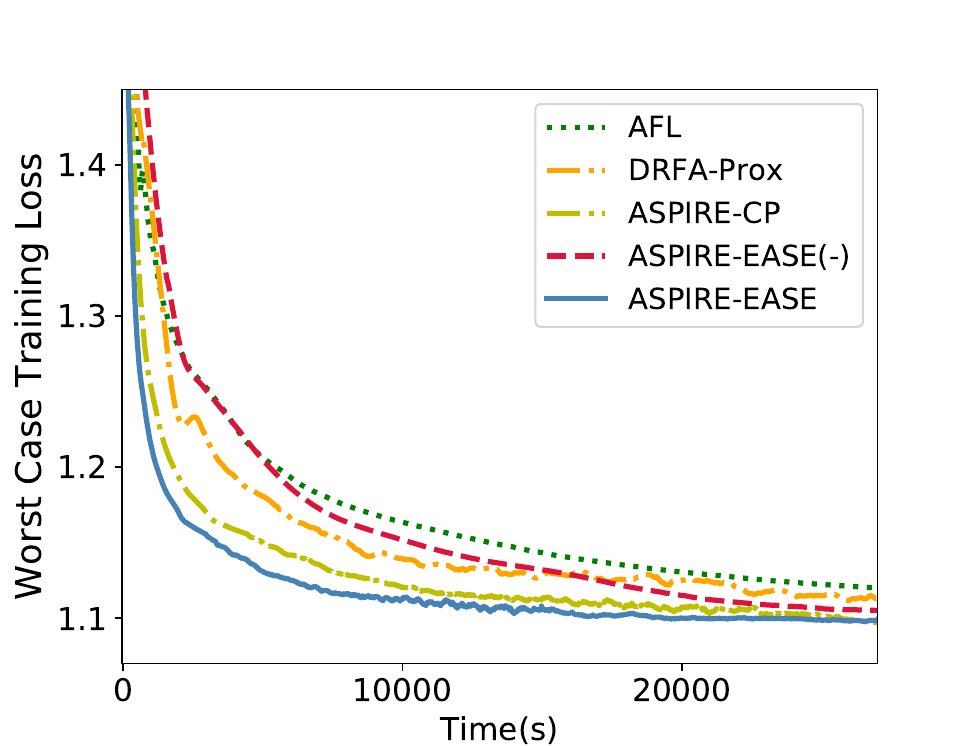}  \vspace{0.01mm}
\end{minipage}}

\caption{Comparison of the convergence time on worst case worker on  (a) Person Activity, (b) SC-MA  datasets.} 
\label{fig:time efficientcy1}
\end{figure}

Within ASPIRE-EASE, the level of robustness can be controlled by adjusting $\Gamma$. Specially, when $\Gamma=0$, we obtain a nominal optimization problem in which no adversarial distribution is considered. The size of the uncertainty set will increase with $\Gamma$ (when $\Gamma \le N$), which enhances the adversarial robustness of the model. As shown in Figures \ref{fig:gamma1} and \ref{fig:gamma2}, the robustness of ASPIRE-EASE can be gradually enhanced when  $\Gamma$ increases.

\begin{figure}[t]
\subfigure[SHL]
{\begin{minipage}{4.2cm}
    \includegraphics[scale=0.285]{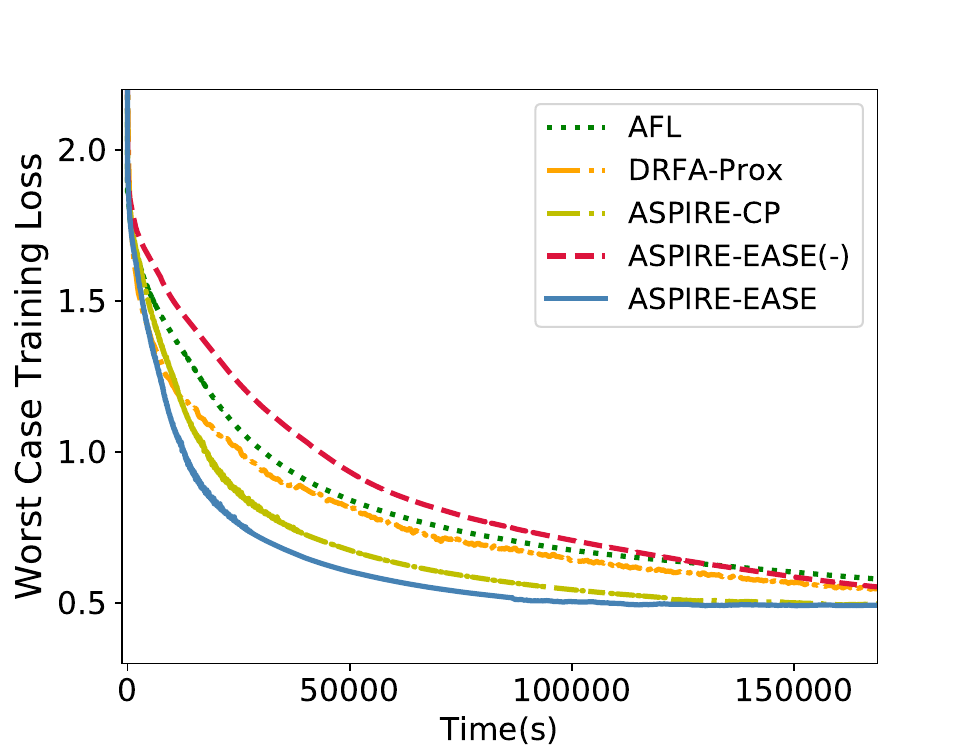}  \vspace{0.01mm}
\end{minipage}}
\subfigure[Fashion MNIST] 
{\begin{minipage}{4.2cm}
    \includegraphics[scale=0.285]{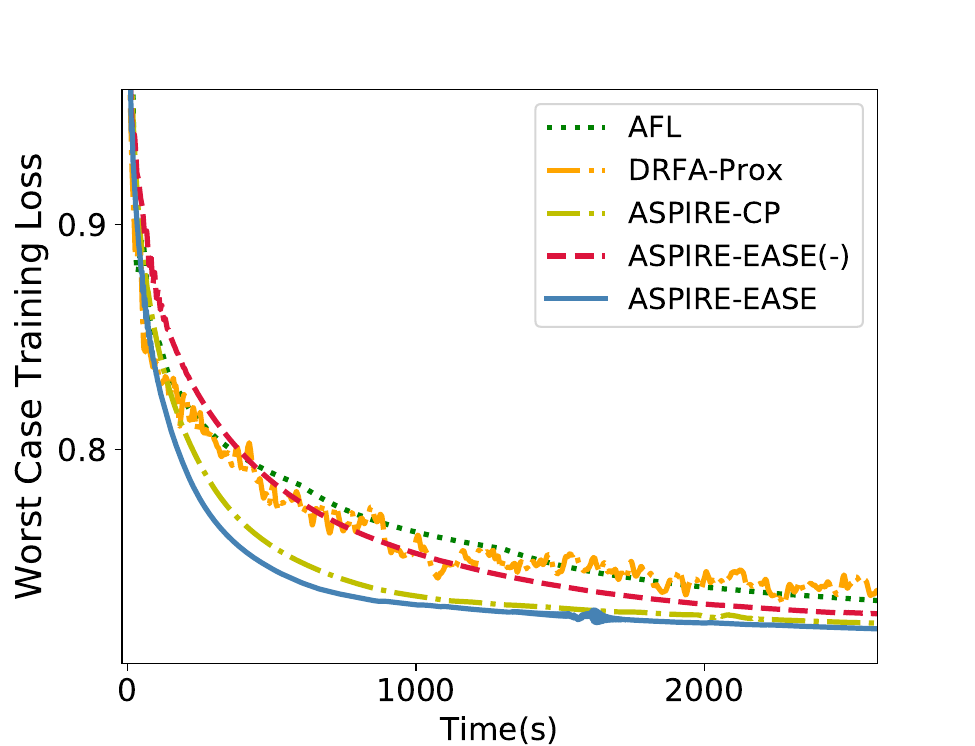}  \vspace{0.01mm}
\end{minipage}}

\caption{Comparison of the convergence time on worst case worker on (a) SHL, (b) Fashion MNIST  datasets.} 
\label{fig:time efficientcy2}
\end{figure}

\begin{figure}[t]
\subfigure[Person Activity]
{\begin{minipage}{4.2cm}
    \includegraphics[scale=0.28]{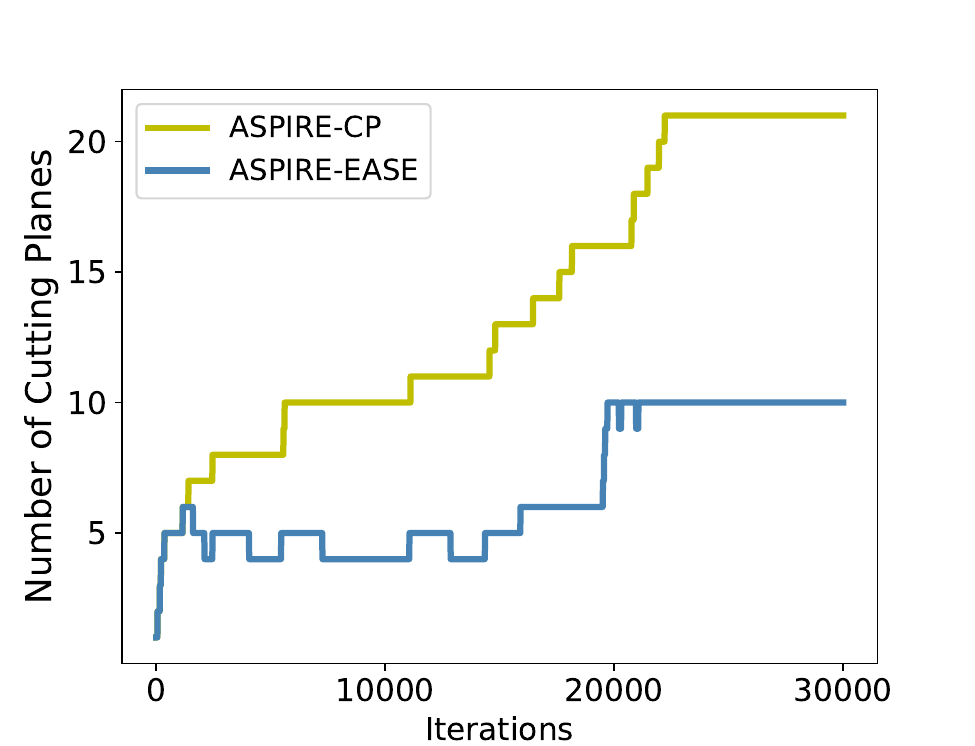}  \vspace{0.01mm}
\end{minipage}}
\subfigure[SC-MA] 
{\begin{minipage}{4.2cm}
    \includegraphics[scale=0.28]{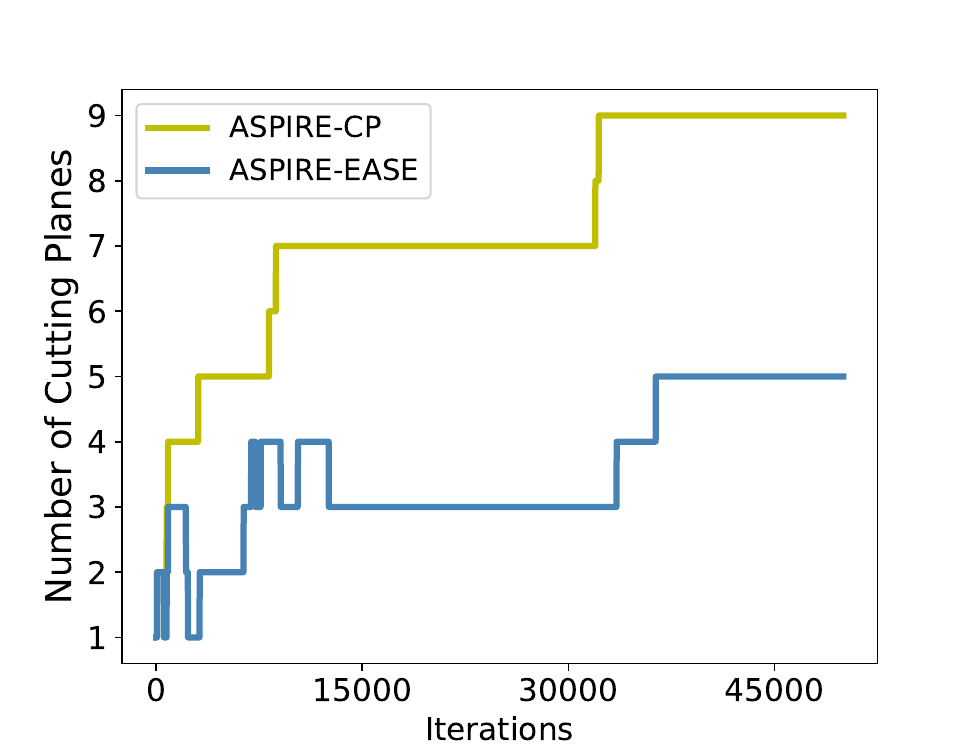}  \vspace{0.01mm}
\end{minipage}}

\caption{Comparison of ASPIRE-CP and ASPIRE-EASE regarding the number of cutting planes on (a) Person Activity, (b) SC-MA  datasets.} 
\label{fig:cutting plane1}
\end{figure}

\begin{figure}[t]
\subfigure[SHL]
{\begin{minipage}{4.2cm}
    \includegraphics[scale=0.28]{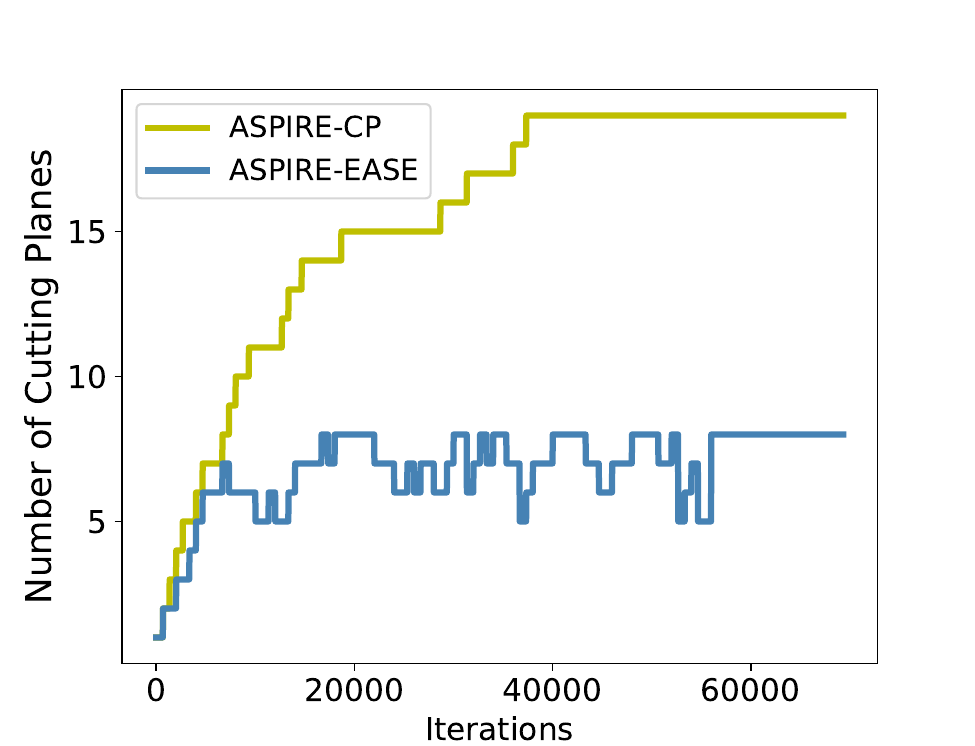}  \vspace{0.01mm}
\end{minipage}}
\subfigure[Fashion MNIST] 
{\begin{minipage}{4.2cm}
    \includegraphics[scale=0.28]{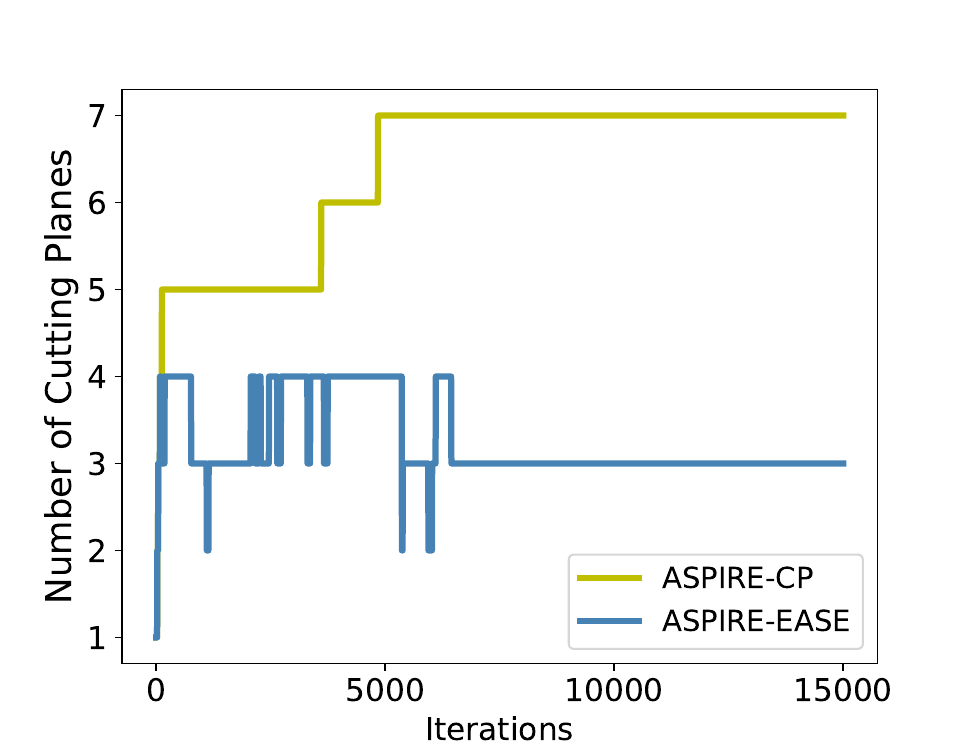}  \vspace{0.01mm}
\end{minipage}}

\caption{Comparison of ASPIRE-CP and ASPIRE-EASE regarding the number of cutting planes on (a) SHL, (b) Fashion MNIST  datasets.} 
\label{fig:cutting plane2}
\end{figure}

\begin{figure}[t]
\subfigure[Test loss vs time]
{\begin{minipage}{4.2cm}
    \includegraphics[scale=0.26]{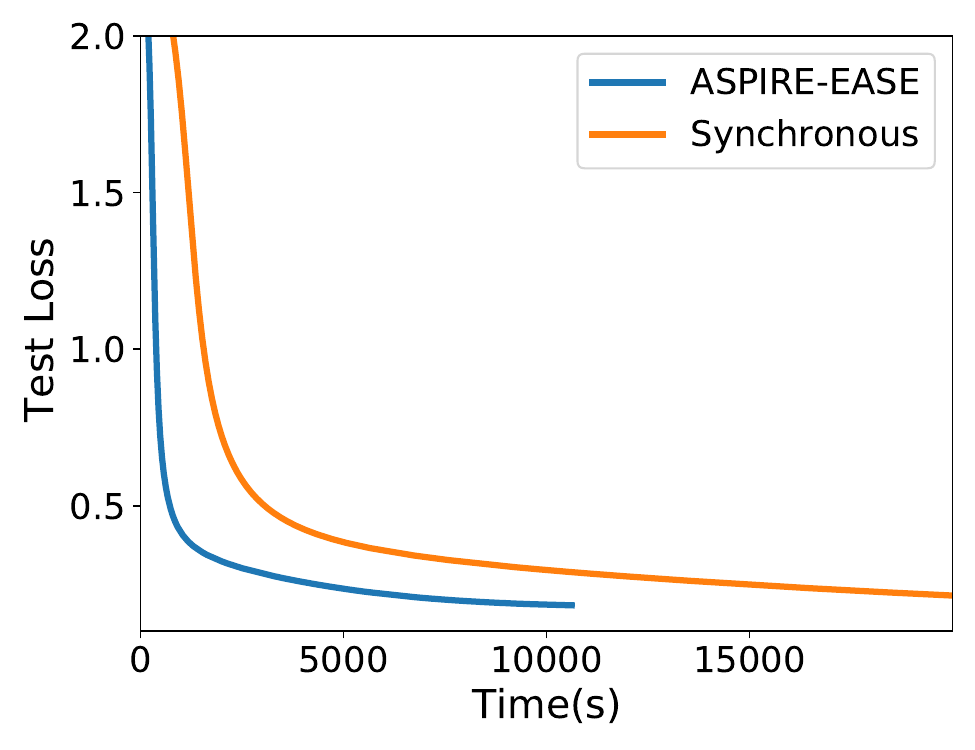}  
\end{minipage}}
\subfigure[Test accuracy vs time] 
{\begin{minipage}{4.2cm}
    \includegraphics[scale=0.26]{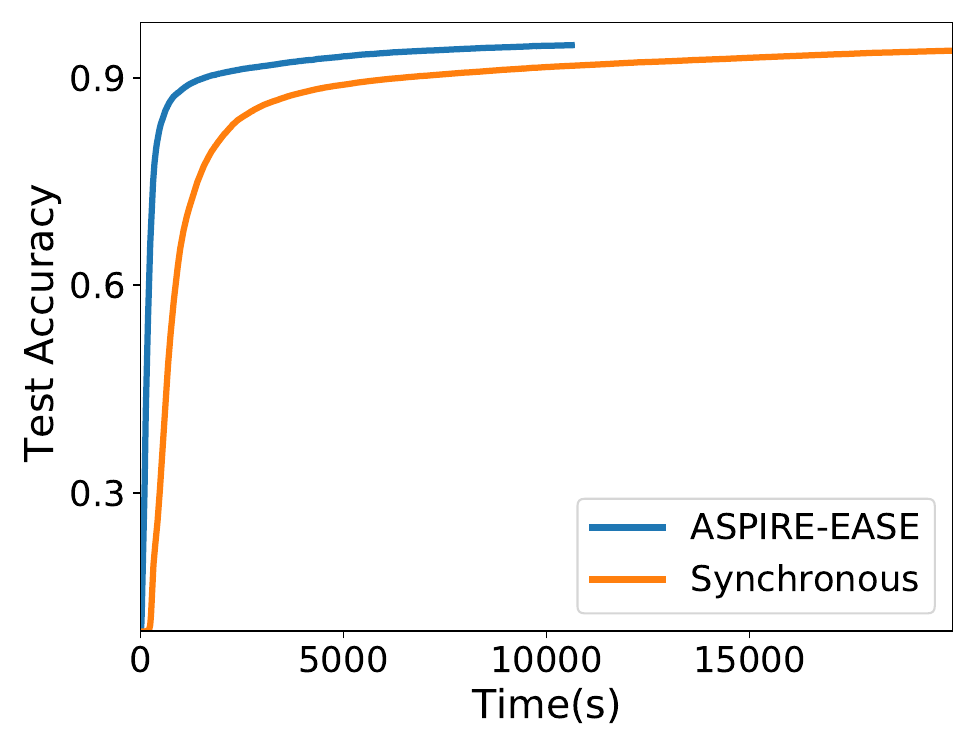}  
\end{minipage}}

\caption{\textcolor{black}{Comparison of ASPIRE-EASE and the synchronous distributed algorithm regarding (a) test loss vs time, and (b) test accuracy vs time on QMNIST dataset with 1000 workers.}} 
\label{fig:sca}
\end{figure}

\subsubsection{Robustness against Malicious Attacks}
To assess the model robustness against malicious attacks, malicious workers with backdoor attacks \citep{bagdasaryan2020backdoor}, which attempt to mislead the model training process, are added to the distributed system.  Following \citep{dai2019backdoor}, we report the success attack rate of backdoor attacks for comparison. It can be calculated by checking how many instances in the backdoor dataset can be misled and categorized into the target labels. Lower success attack rates indicate more robustness against backdoor attacks.  The comparison results are summarized in Table~\ref{tab:SAR} and more detailed settings of backdoor attacks are available in Appendix C in Supplementary Materials.  In Table~\ref{tab:SAR}, we observe that AFL can be attacked easily since it could assign higher weights to malicious workers. Compared to AFL, FedAvg and ${\rm{Mix}}{\rm{_{Even}}}$ achieve relatively lower success attack rates since they assign equal weights to the malicious workers and other workers. DRFA-Prox can achieve even lower success attack rates since it can leverage the prior distribution to assign lower weights for malicious workers.  The proposed ASPIRE-EASE achieves the lowest success attack rates since it can leverage the prior distribution more effectively. Specifically, it will assign lower weights to malicious workers with tight theoretical guarantees.

\renewcommand\tabcolsep{5pt}
\renewcommand\arraystretch{1.4}
\begin{table}[t]
\caption{Performance comparisons about the success attack rate ($\%$) $\downarrow$. The boldfaced digits represent the best results.}
\label{tab:SAR}
\centering
\scalebox{0.9}{
\begin{tabular}{l|cccc}
\toprule
Model    & SHL     & Person Activity         & SC-MA     & Fashion MNIST        \\ \hline
${\rm{Mix}}{\rm{_{Even}}}$ & 36.21±2.23 &  34.32±2.18  & 52.14±2.89   & 83.18±2.07  \\ 
FedAvg \citep{mcmahan2017communication} & 38.15±3.02  & 33.25±2.49 & 55.39±3.13    &  82.04±1.84\\ 
AFL \citep{mohri2019agnostic} & 68.63±4.24  &  43.66±3.87 & 75.81±4.03  &  90.04±2.52 \\
DRFA-Prox \citep{deng2021distributionally} & 21.23±3.63  &  27.27±3.31 & 30.79±3.65  &  63.24±2.47  \\
\hline
ASPIRE-EASE & \textbf{9.17±1.65}  &  \textbf{22.36±2.33}  & \textbf{14.51±3.21}  & \textbf{45.10±1.64}
\cr
\bottomrule  
\end{tabular}}
\end{table}

\subsubsection{Efficiency}
In Figures \ref{fig:time efficientcy1} and \ref{fig:time efficientcy2}, we compare the efficiency of the proposed ASPIRE-EASE with AFL, DRFA-Prox, and two variants of ASPIRE-EASE, namely: ASPIRE-CP (ASPIRE with the cutting plane method) and ASPIRE-EASE(-) (ASPIRE-EASE without the asynchronous setting). Based on the comparison, we can observe that the proposed ASPIRE-EASE generally achieves higher efficiency than baseline methods and its two variants. This is because 1) compared with AFL, DRFA-Prox, and ASPIRE-EASE(-), ASPIRE-EASE is an \textit{asynchronous} algorithm in which the master updates its parameters only after receiving the updates from active workers instead of all workers; 2) the proposed ASPIRE-EASE is a \textit{single-loop} algorithm with \textit{efficient} computations at each iteration, unlike DRFA-Prox, the master in ASPIRE-EASE only needs to communicate with active workers once per iteration; 3) compared with ASPIRE-CP, ASPIRE-EASE utilizes \textit{active set method} instead of cutting plane method, which is more efficient. It is seen from Figures \ref{fig:time efficientcy1} and \ref{fig:time efficientcy2} that, the convergence speed of ASPIRE-EASE mainly benefits from the asynchronous setting.

\textcolor{black}{\subsubsection{Scalability on Large-scale Problem}Distributed optimization serves as an effective framework for addressing large-scale optimization problems, as it enables the decomposition of complex problems into smaller subproblems that can be solved in parallel. Moreover, compared to synchronous distributed algorithms, asynchronous algorithms are better suited for large-scale environments, as they are not affected by the straggler problem and exhibit greater robustness to device failures. In this work, we propose an efficient single-loop asynchronous distributed algorithm that is well-suited for large-scale optimization tasks. To empirically evaluate its scalability, we conduct experiments on a large-scale distributed system involving 1,000 workers. In the experiment, the QMNIST dataset is partitioned into 1,000 sub-datasets, each assigned to a different worker, and 5 workers are set as stragglers. As shown in Figure \ref{fig:sca}, the proposed asynchronous algorithm effectively addresses the large-scale learning problem.}

\subsubsection{Ablation Study}
For ASPIRE, compared with cutting plane method, EASE is more efficient since it considers removing the inactive cutting planes. To demonstrate the efficiency of EASE, we first compare ASPIRE-EASE with ASPIRE-CP concerning the number of cutting planes used during the training. In Figures \ref{fig:cutting plane1} and \ref{fig:cutting plane2}, we can observe that ASPIRE-EASE uses fewer cutting planes than ASPIRE-CP, thus is more efficient. The convergence speed of ASPIRE-EASE and ASPIRE-CP in Figures \ref{fig:time efficientcy1} and \ref{fig:time efficientcy2} also suggests that ASPIRE-EASE converges much faster than ASPIRE-CP.

\textcolor{black}{
\section{Limitations and Future Work}
In this work, we focus on the federated distributionally robust optimization (FDRO). While the federated learning framework avoids the need to upload data to a central server, which offers a degree of privacy protection, the level of privacy remains limited \cite{geiping2020inverting}. Therefore, it is critical to investigate privacy-preserving FDRO in future work. Inspired \cite{ma2025cellular}, where the federated learning with differential privacy is formulated as a min-max optimization problem (which can also be viewed as a bilevel optimization problem), we extend this perspective and formulate privacy-preserving FDRO as a federated trilevel optimization problem. In comparison to the bilevel FDRO studied in this work, the trilevel privacy-preserving FDRO requires more computational costs due to its more intricate hierarchical structure, as discussed in \cite{jiao2024provably,jiao2024unlocking}. This represents a trade-off between the strength of privacy guarantees and computational efficiency.}

\textcolor{black}{The federated distributionally robust optimization under the setting that all participating workers/clients are reliable and willing to engage in training is considered in this work. Recently, increasing attention has been paid to scenarios where clients may be reluctant to participate or share their model updates without sufficient compensation \cite{zhan2021survey,ding2025incentivized,pandey2019incentivize}. In such settings, the design of effective incentivization mechanisms becomes crucial. For instance, \cite{siraj2022incentives} proposes an incentive mechanism based on the theory of Colonel Blotto games, while \cite{zhao2021incentive} adopts the framework of Stackelberg game to address the same challenge. In future work, we plan to design effective incentivization mechanisms to address federated distributionally robust optimization in scenarios where client participation is not guaranteed.
}

\section{Conclusion} \label{conclusion}

In this paper, we introduce the ASPIRE-EASE method for efficiently addressing the federated distributionally robust optimization (FDRO) problems characterized by non-convex objectives. Additionally, $CD$-norm uncertainty set has been proposed to effectively integrate the prior distribution into the problem formulation, thereby enabling the flexible modulation of the robustness of DRO. To further enhance the proposed framework, we integrate various uncertainty sets and conduct a comprehensive theoretical analysis, delving into the computational complexities associated with each uncertainty set. To enhance the efficiency of the proposed method, ASPIRE-ADP, which is a methodology capable of dynamically adjusting the number of active workers, is proposed. Theoretical analyses have been systematically performed to evaluate the convergence properties, iteration complexity, and communication complexity of the proposed method. Extensive experimental results illustrate that the proposed method exhibits noteworthy empirical efficacy across diverse real-world datasets, establishing its effectiveness in addressing the DRO  problems through a fully distributed and asynchronous manner.

 
\ifCLASSOPTIONcaptionsoff
  \newpage
\fi

\bibliographystyle{ieeetr}
\bibliography{ref}

\end{document}